\documentclass[11pt]{article}

\usepackage{latexsym}
\usepackage{amssymb}
\usepackage{amsthm}
\usepackage{amscd}
\usepackage{amsmath} 
\usepackage[all]{xy}
\usepackage{graphicx}

\newtheorem{theorem}{Theorem}[section]

\numberwithin{equation}{section}

\xyoption{dvips}

\CompileMatrices

\def\F{\mathbb{F}}   \def\C{\mathbb{C}} 
   
 \def\dim{\mathrm{dim}} \def\Hom{\mathrm{Hom}} 
\def\End{\mathrm{End}} \def\Z{\mathbb{Z}} 
  \def\H{\mathcal{H}} \def\Ind{\mathrm{Ind}} \def\GL{\mathrm{GL}}   \def\Res{\mathrm{Res}} \def\dd{\displaystyle} \def\Inf{\mathrm{Inf}} \def\ch{\mathrm{ch}} \def\spanning{\textnormal{-span}}
\def\vphi{\varphi} \def\U{\mathcal{U}} \def\B{\mathcal{B}}
  \def\wt{\mathrm{wt}}\def\Indf{\mathrm{Indf}}
 \def\L{\mathcal{L}} \def\sh{\mathrm{sh}}\def\hgt{\mathrm{ht}}
 \def\Card{\mathrm{Card}}  \def\ss{\scriptstyle}
\def\mcP{\mathcal{P}}

\makeatletter
\renewcommand{\@makefnmark}{\mbox{\textsuperscript{}}}
\makeatother

\title{Unipotent Hecke algebras of $\GL_n(\F_q)$}
\author{Nathaniel Thiem\\University of Wisconsin-Madison\\ 480 Lincoln Drive \\ Madison, WI 53706\\ {\small\textsf{thiem@math.wisc.edu}}}
\date{ }

\begin{document}

\maketitle

\begin{abstract} This paper describes a family of Hecke algebras $\H_\mu=\End_G(\Ind_U^G(\psi_\mu))$, where $U$ is the subgroup of unipotent upper-triangular matrices of $G=\GL_n(\F_q)$ and $\psi_\mu$ is a linear character of $U$.  The main results combinatorially index a basis of $\H_\mu$, provide a large commutative subalgebra of $\H_\mu$,  and after describing the combinatorics associated with the representation theory of $\H_\mu$, generalize the RSK correspondence that is typically found in the representation theory of the symmetric group.  
\end{abstract}

\begin{section}{Introduction\protect\protect}

Iwahori \cite{Iw} and Iwahori-Matsumoto \cite{IM} introduced the Iwahori-Hecke algebra as a first step in classifying the irreducible representations of finite Chevalley groups and reductive $p$-adic Lie groups.  Subsequent work (e.g. \cite{CuS} \cite{KL} \cite{LV}) has established Hecke algebras as fundamental tools in the representation theory of Lie groups and Lie algebras, and advances on subfactors and quantum groups by Jones \cite{Jo1}, Jimbo \cite{Ji}, and Drinfeld \cite{Dr} gave Hecke algebras a central role in knot theory \cite{Jo2}, statistical mechanics \cite{Jo3}, mathematical physics, and operator algebras.  This paper considers a generalization of the classical Iwahori-Hecke algebra obtained by replacing the Borel subgroup $B$ with the unipotent subgroup $U$.

A \emph{Hecke algebra} $\H=\H(G,U,M)$ is the centralizer algebra
$$\H=\End_G(\Ind_U^G(M)),$$
where $G$ is a finite group, $U$ is a subgroup of $G$ and $M$ is a simple $U$-module.  This paper addresses the cases where
$$G=\GL_n(\F_q),\qquad U=\left\{\left(\begin{array}{ccc} 1 & & * \\ & \ddots & \\ 0 & & 1\end{array}\right)\right\},\qquad \textnormal{and}\qquad \dim(M)=1.$$
For each $M$ there exist a partition $\mu=(\mu_1,\mu_2,\ldots,\mu_\ell)$ and a linear character $\psi_\mu:U\rightarrow \C^*$ such that
$$  \End_G(\Ind_U^G(M))\cong\End_G(\Ind_U^G(\psi_\mu))=\H_\mu.$$
In studying these algebras, this paper combines generalizations of Iwahori-Hecke algebra techniques with combinatorial tools related to the representation theory of the symmetric group (e.g. partitions, tableaux, and the RSK correspondence).

The main results of this paper are
\begin{enumerate}
\item[(a)]  Let $\{T_v\ :\ v\in N_\mu\}$ be the standard double-coset basis for $\H_\mu$.  Then there is a bijection
$$N_\mu\longleftrightarrow M_\mu=\left\{\begin{array}{p{6.5cm}} $\ell\times \ell$ matrices with monic polynomial entries $a_{ij}(X)\in \F_q[X]$ such that $a_{ij}(0)\neq 0$ and both degree row sums and degree column sums are equal to $\mu$\end{array}\right\}$$
(See Section 3 for details).
\item[(b)]  Let the set $\hat\H_\mu$ index the irreducible $\H_\mu$-modules $\H_\mu^\lambda$ and the set $\hat\H_\mu^\lambda$ index a basis of the irreducible module $\H_\mu^\lambda$.  There is a combinatorial bijection 
$$N_\mu \longleftrightarrow \left\{\begin{array}{c}\textnormal{Pairs } (P,Q) \textnormal{ such}\\ \textnormal{that } P,Q\in \H_\mu^\lambda, \lambda\in \hat\H_\mu\end{array}\right\}$$
that generalizes the classical RSK correspondence and gives a combinatorial realization of the representation theoretic identity
$$\dim(\H_\mu)=\sum_{\lambda\in \hat\H_\mu} \dim(\H_\mu^\lambda)^2.$$
\item[(c)] The algebra $\H_\mu$ has a large commutative subalgebra
$$\L_\mu\cong \H_{(\mu_1)}\otimes \H_{(\mu_2)}\otimes\cdots \otimes \H_{(\mu_\ell)},$$ 
whose presence suggests a weight space decomposition of $\H_\mu$-modules.

\end{enumerate}

Section 2 reviews some basic results used in this paper, including Hecke algebras, partitions and the classical RSK correspondence.  Section 3 defines $\H_\mu$ and gives an explicit construction of the map $N_\mu\leftrightarrow M_\mu$.  After using Zelevinsky's theorem \cite[Theorem 12.1]{Z} to give a combinatorial description of the sets $\hat\H_\mu$ and $\hat\H_\mu^\lambda$, Section 4 proves the bijection in (b).  Section 5 gives a proof of Zelevinsky's theorem.  This paper concludes in Section 6 by providing the subalgebra $\L_\mu\subseteq \H_\mu$ and describing a corresponding weight space decomposition of $\H_\mu$-modules.   

Both the Yokonuma algebra $\H_{(1^n)}$ \cite{Y2} and the Hecke algebra $\H_{(n)}$ associated to the Gelfand-Graev representation of $G$ \cite{St} are examples of unipotent Hecke algebras.  In \cite{Y2}, Yokonuma gave a presentation of $\H_{(1^n)}$ that generalized the usual presentation of the classical Iwahori-Hecke algebra, and recently Jujumaya \cite{Ju} constructed an alternate set of generators and relations.  However, a presentation for arbitrary $\H_\mu$ is still unknown.  Even the commutative algebra $\H_{(n)}$ \cite{St,Y1} does not yet have a ``nice" set of multiplication relations.

The representation theory of the Gelfand-Graev Hecke algebra is closely related to Green polynomials \cite{Cu} and, in the $\GL_2$ case, to Kloosterman sums \cite{CS}.  On the other hand, the representation theory of the Yokonuma algebra generalizes that of the classical Iwahori-Hecke algebra.  In what promises to be a combinatorially rich area of study, analyzing the combinatorial representation theory of the $\H_\mu$ and their general type analogues should have an impact similar in scope to the applications of the classical Iwahori-Hecke algebra.

\vspace{.25cm}

\noindent\textbf{Acknowledgements.}  This paper will be a portion of my Ph.D. thesis.  In developing these results, I have enjoyed the supportive environment of the University of Wisconsin-Madison mathematics department, the time supplied by several grants (VIGRE DMS-9819788, NSF DMS-0097977, and NSA MDA904-01-1-0032), and above all the patient help and insights of my advisor Arun Ram.

\end{section}

\begin{section}{Preliminaries: Hecke algebras, some combinatorics of the symmetric group, and $\GL_n(\F_q)$}

\noindent\textbf{Hecke algebras.}  Let $U$ be a subgroup of a finite group $G$.  If $M$ is an irreducible $U$-module, then the \emph{Hecke algebra} $\H=\H(G,U,M)$ is
$$\H=\End_G\left(\Ind_U^G(M)\right)\cong e \C G e, $$
where $e$ is an idempotent of $\C U$ such that $M\cong \C U e$ \cite[(3.19)]{CR}.  If $N$ is a set of double coset representatives for the cosets $U\backslash G/U$, then the set
\begin{equation}\label{Tbasis}\{T_v=eve\ :\ v\in N, eve\neq 0\}\end{equation}
is a basis for $\H$ \cite[(11.30)]{CR}.

Let $\hat{\H}$ be an indexing set for the irreducible $\H$-modules $\H^\lambda$.  As a $(G,\H)$-bimodule, 
$$\C G e\cong \Ind_U^G(M)\cong \bigoplus_{\lambda\in \hat\H} G^\lambda\otimes \H^\lambda,$$
where the $G^\lambda$ are the irreducible constituents of $\Ind_U^G(M)$ \cite[Thm 3.3.7]{GW}; it follows that 
\begin{equation}\label{SW} \dim(\H^\lambda)=\textnormal{ multiplicity of } G^\lambda\textnormal{ in the } G\text{-module } \Ind_U^G(M). \end{equation}

\vspace{.5cm}

\noindent\textbf{Compositions, partitions and tableaux.}  A \emph{composition } $\mu=(\mu_1,\mu_2,\ldots,\mu_r)$ is a sequence of positive integers.  The \emph{size} of $\mu$ is $\vert \mu\vert=\mu_1 + \mu_2 + \cdots + \mu_r$, the \emph{length} of $\mu$ is $\ell(\mu)=r$ and
\begin{equation}\label{Bmu}\B_\mu=\{\mu_1,\mu_1+\mu_2,\ldots,\mu_1+\mu_2+\cdots +\mu_r\}.\end{equation}
If $|\mu|=n$, then $\mu$ is \emph{a composition of $n$} and we write $\mu\models n$.  View $\mu$ as a collection of boxes aligned to the left; for example, if
$$\mu=(2,5,3,4)=\xy<0cm,.8cm>\xymatrix@R=.4cm@C=.4cm{ 
*={} & *={} \ar @{-} [l] & *={} \ar @{-} [l]\\
*={}\ar @{-} [u] & *={} \ar @{-} [l] \ar @{-} [u] & *={} \ar @{-} [l] \ar @{-} [u] & *={} \ar @{-} [l] & *={} \ar @{-} [l] & *={} \ar @{-} [l] \\
*={}\ar @{-} [u] & *={} \ar @{-} [l] \ar @{-} [u] & *={} \ar @{-} [l] \ar @{-} [u] & *={}\ar @{-} [u] \ar @{-} [l] & *={} \ar @{-} [l] \ar @{-} [u] & *={} \ar @{-} [l] \ar @{-} [u]  \\
*={}\ar @{-} [u] & *={} \ar @{-} [l] \ar @{-} [u] & *={} \ar @{-} [l] \ar @{-} [u]  & *={} \ar @{-} [u] \ar @{-} [l] & *={} \ar @{-} [l]\\
*={}\ar @{-} [u] & *={} \ar @{-} [l] \ar @{-} [u] & *={} \ar @{-} [l] \ar @{-} [u] & *={} \ar @{-} [u] \ar @{-} [l] & *={} \ar @{-} [l] \ar @{-} [u]} \endxy ,$$
then $|\mu |=14$, $\ell(\mu)=4$ and $\B_\mu=\{2,7,10,14\}$. Alternatively, $\B_\mu$ coincides with the numbers in the boxes at the end of the rows in the diagram
$$\xy<0cm,.8cm>\xymatrix@R=.4cm@C=.4cm{ 
*={} & *={} \ar @{-} [l] & *={} \ar @{-} [l] \ar @/^1pt/ @{-} [l] \ar @/_1pt/ @{-} [d]\\
*={}\ar @{-} [u] & *={} \ar @{-} [l] \ar @{-} [u] \ar @{} [ul]|{1} \ar @/^1pt/ @{-} [r] \ar @/_1pt/ @{-} [u] & *={} \ar @{-} [l] \ar @{-} [u] \ar @{} [ul]|{2} & *={} \ar @{-} [l] & *={} \ar @{-} [l] & *={} \ar @{-} [l] \ar @/^1pt/ @{-} [l] \ar @/_1pt/ @{-} [d] \\
*={}\ar @{-} [u] & *={} \ar @{-} [l] \ar @{-} [u] \ar @{} [ul]|{3} & *={} \ar @{-} [l] \ar @{-} [u] \ar @{} [ul]|{4} & *={}\ar @{-} [u] \ar @{-} [l] \ar @{} [ul]|{5} \ar @/^1pt/ @{-} [l] \ar @/_1pt/ @{-} [d] & *={} \ar @{-} [l] \ar @{-} [u] \ar @{} [ul]|{6} \ar @/^1pt/ @{-} [r] \ar @/_1pt/ @{-} [u] & *={} \ar @{-} [l] \ar @{-} [u] \ar @{} [ul]|{7} \\
*={}\ar @{-} [u] & *={} \ar @{-} [l] \ar @{-} [u] \ar @{} [ul]|{8} & *={} \ar @{-} [l] \ar @{-} [u]  \ar @{} [ul]|{9} \ar @/^1pt/ @{-} [r] \ar @/_1pt/ @{-} [u] & *={} \ar @{-} [u] \ar @{-} [l] \ar @{} [ul]|{10} & *={} \ar @{-} [l] \ar @/^1pt/ @{-} [l] \ar @/_1pt/ @{-} [d] \\
*={}\ar @{-} [u] & *={} \ar @{-} [l] \ar @{-} [u] \ar @{} [ul]|{11} & *={} \ar @{-} [l] \ar @{-} [u] \ar @{} [ul]|{12} & *={} \ar @{-} [u] \ar @{-} [l] \ar @{} [ul]|{13} \ar @/^1pt/ @{-} [r] \ar @/_1pt/ @{-} [u] & *={} \ar @{-} [l] \ar @{-} [u] \ar @{} [ul]|{14}} \endxy\ .$$

A \emph{partition} $\nu=(\nu_1,\nu_2,\ldots,\nu_r)$ is a composition where $\nu_1\geq \nu_2\geq \cdots \geq \nu_r>0.$  If $|\nu |=n$, then $\nu$ is \emph{a partition of $n$} and we write $\nu\vdash n$.  Let
\begin{equation}\label{partitions} \mcP=\{\text{partitions}\}\qquad \text{and} \qquad \mcP_n=\{\nu\vdash n\}.\end{equation}
Suppose $\nu\in\mcP$.  The \emph{conjugate partition} $\nu^\prime=(\nu_1^\prime,\nu_2^\prime,\ldots,\nu_\ell^\prime)$ is given by
$$\nu_i^\prime=\Card\{j\ :\ \nu_j\geq i\}.$$
In terms of diagrams, $\nu^\prime$ is the collection of boxes obtained by flipping $\nu$ across its main diagonal.  For example,
$$\textnormal{if}\qquad \nu=
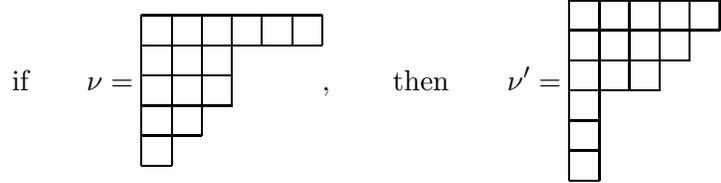
\xy<0cm,.9cm>\xymatrix@R=.4cm@C=.4cm{
*={} & *={} \ar @{-} [l] & *={} \ar @{-} [l] & *={} \ar @{-} [l] & *={} \ar @{-} [l] & *={} \ar @{-} [l]  & *={} \ar @{-} [l] \\
*={} \ar @{-} [u] &   *={} \ar @{-} [l] \ar@{-} [u] &   *={} \ar @{-} [l] \ar@{-} [u] &   *={} \ar @{-} [l] \ar@{-} [u] &   *={} \ar @{-} [l] \ar@{-} [u] &   *={} \ar @{-} [l] \ar@{-} [u] & *={}\ar @{-} [l] \ar @{-} [u]\\
*={} \ar @{-} [u] &   *={} \ar @{-} [l] \ar@{-} [u] &   *={} \ar @{-} [l] \ar@{-} [u] &   *={} \ar @{-} [l] \ar@{-} [u] \\
*={} \ar @{-} [u] &   *={} \ar @{-} [l] \ar@{-} [u] &   *={} \ar @{-} [l] \ar@{-} [u] &   *={} \ar @{-} [l] \ar@{-} [u] \\
*={} \ar @{-} [u] &   *={} \ar @{-} [l] \ar@{-} [u] &   *={} \ar @{-} [l] \ar@{-} [u] \\
*={} \ar @{-} [u] &   *={} \ar @{-} [l] \ar@{-} [u] }\endxy,\qquad \textnormal{then}  \qquad
\nu^\prime=
\xy<0cm,1.1cm>\xymatrix@R=.4cm@C=.4cm{
*={} & *={} \ar @{-} [l] & *={} \ar @{-} [l] & *={} \ar @{-} [l] & *={} \ar @{-} [l] & *={} \ar @{-} [l]  \\
*={} \ar @{-} [u] &   *={} \ar @{-} [l] \ar@{-} [u] &   *={} \ar @{-} [l] \ar@{-} [u] &   *={} \ar @{-} [l] \ar@{-} [u] &   *={} \ar @{-} [l] \ar@{-} [u] & *={} \ar @{-} [l]  \ar @{-} [u] \\
*={} \ar @{-} [u] &   *={} \ar @{-} [l] \ar@{-} [u] &   *={} \ar @{-} [l] \ar@{-} [u] &   *={} \ar @{-} [l] \ar@{-} [u] & *={} \ar @{-} [u] \ar @{-} [l]\\
*={} \ar @{-} [u] &   *={} \ar @{-} [l] \ar@{-} [u] &   *={} \ar @{-} [l] \ar@{-} [u] &   *={} \ar @{-} [l] \ar@{-} [u] \\
*={} \ar @{-} [u] &   *={} \ar @{-} [l] \ar@{-} [u] \\
*={} \ar @{-} [u] &   *={} \ar @{-} [l] \ar@{-} [u] \\
*={} \ar @{-} [u] &   *={} \ar @{-} [l] \ar@{-} [u]  }\endxy .$$

A \emph{column strict tableau $Q$ of shape $\nu$} is a filling of the boxes of $\nu$ by positive integers such that 
\begin{enumerate}
\item[(a)] the entries strictly increase along columns,
\item[(b)] the entries weakly increase along rows.
\end{enumerate}
The \emph{weight of $Q$} is the composition $\wt(Q)=(\wt(Q)_1,\wt(Q)_2,\ldots )$ given by
$$\wt(Q)_i=\textnormal{number of } i \textnormal{ in } Q.$$ 
For example,
$$ Q=
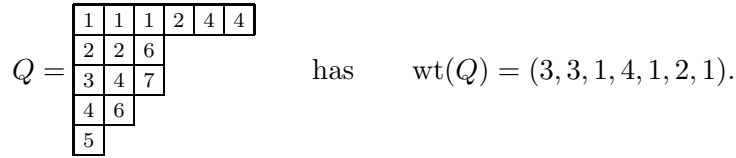
\xy<0cm,.9cm>\xymatrix@R=.4cm@C=.4cm{
*={} & *={} \ar @{-} [l] & *={} \ar @{-} [l] & *={} \ar @{-} [l] & *={} \ar @{-} [l] & *={} \ar @{-} [l]  & *={} \ar @{-} [l] \\
*={} \ar @{-} [u] &   *={} \ar @{-} [l] \ar@{-} [u] \ar @{} [ul]|{1} &   *={} \ar @{-} [l] \ar@{-} [u] \ar @{} [ul]|{1} &   *={} \ar @{-} [l] \ar@{-} [u] \ar @{} [ul]|{1} &   *={} \ar @{-} [l] \ar@{-} [u] \ar @{} [ul]|{2} &   *={} \ar @{-} [l] \ar@{-} [u] \ar @{} [ul]|{4} & *={}\ar @{-} [l] \ar @{-} [u] \ar @{} [ul]|{4}\\
*={} \ar @{-} [u] &   *={} \ar @{-} [l] \ar@{-} [u] \ar @{} [ul]|{2} &   *={} \ar @{-} [l] \ar@{-} [u]  \ar @{} [ul]|{2} &   *={} \ar @{-} [l] \ar@{-} [u] \ar @{} [ul]|{6} \\
*={} \ar @{-} [u]  &   *={} \ar @{-} [l] \ar@{-} [u] \ar @{} [ul]|{3} &   *={} \ar @{-} [l] \ar@{-} [u] \ar @{} [ul]|{4} &   *={} \ar @{-} [l] \ar@{-} [u] \ar @{} [ul]|{7} \\
*={} \ar @{-} [u] &   *={} \ar @{-} [l] \ar@{-} [u] \ar @{} [ul]|{4} &   *={} \ar @{-} [l] \ar@{-} [u] \ar @{} [ul]|{6}\\
*={} \ar @{-} [u] &   *={} \ar @{-} [l] \ar@{-} [u] \ar @{} [ul]|{5}}\endxy\qquad \textnormal{has} \qquad  \wt(Q)=(3,3,1,4,1,2,1).$$

Suppose $\nu,\mu$ are partitions.  If $\nu_i\geq \mu_i$ for all $1\leq i\leq \ell(\mu)$, then
the \emph{skew partition} $\nu/\mu$ is given by 
$$\nu/\mu=(\nu_1-\mu_1,\nu_2-\mu_2,\ldots, \nu_{\ell(\nu)}-\mu_{\ell(\nu)}),$$
where $\mu_k=0$ for all $k\geq \ell(\mu)$.  In terms of boxes, represent $\nu/\mu$ by removing $\mu$ from the upper left-hand corner of the diagram $\nu$, so if
$$\nu=
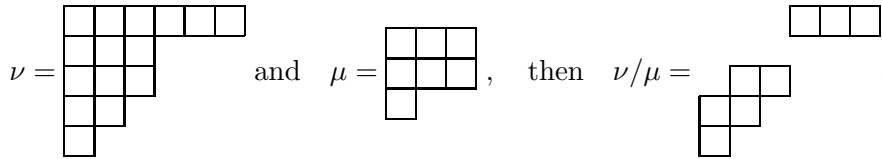
\xy<0cm,.9cm>\xymatrix@R=.4cm@C=.4cm{
*={} & *={} \ar @{-} [l] & *={} \ar @{-} [l] & *={} \ar @{-} [l] & *={} \ar @{-} [l] & *={} \ar @{-} [l]  & *={} \ar @{-} [l] \\
*={} \ar @{-} [u] &   *={} \ar @{-} [l] \ar@{-} [u] &   *={} \ar @{-} [l] \ar@{-} [u] &   *={} \ar @{-} [l] \ar@{-} [u] &   *={} \ar @{-} [l] \ar@{-} [u] &   *={} \ar @{-} [l] \ar@{-} [u] & *={}\ar @{-} [l] \ar @{-} [u]\\
*={} \ar @{-} [u] &   *={} \ar @{-} [l] \ar@{-} [u] &   *={} \ar @{-} [l] \ar@{-} [u] &   *={} \ar @{-} [l] \ar@{-} [u] \\
*={} \ar @{-} [u] &   *={} \ar @{-} [l] \ar@{-} [u] &   *={} \ar @{-} [l] \ar@{-} [u] &   *={} \ar @{-} [l] \ar@{-} [u] \\
*={} \ar @{-} [u] &   *={} \ar @{-} [l] \ar@{-} [u] &   *={} \ar @{-} [l] \ar@{-} [u] \\
*={} \ar @{-} [u] &   *={} \ar @{-} [l] \ar@{-} [u] }\endxy \textnormal{ and} \quad \mu=\xy<0cm,.6cm>\xymatrix@R=.4cm@C=.4cm{
*={} & *={} \ar @{-} [l] & *={} \ar @{-} [l] & *={} \ar @{-} [l]\\
*={} \ar @{-} [u] &   *={} \ar @{-} [l] \ar@{-} [u] &   *={} \ar @{-} [l] \ar@{-} [u] &   *={} \ar @{-} [l] \ar@{-} [u] \\
*={} \ar @{-} [u] &   *={} \ar @{-} [l] \ar@{-} [u] &   *={} \ar @{-} [l] \ar@{-} [u] &   *={} \ar @{-} [l] \ar@{-} [u] \\
*={} \ar @{-} [u] &   *={} \ar @{-} [l] \ar@{-} [u]  }\endxy\ , \quad \textnormal{then} \quad \nu/\mu=\xy<0cm,.9cm>\xymatrix@R=.4cm@C=.4cm{
*={} & *={}  & *={} & *={}  & *={} \ar @{-} [l] & *={} \ar @{-} [l]  & *={} \ar @{-} [l] \\
*={}  &   *={}  &   *={}  &   *={}  \ar@{-} [u] &   *={} \ar @{-} [l] \ar@{-} [u] &   *={} \ar @{-} [l] \ar@{-} [u] & *={}\ar @{-} [l] \ar @{-} [u]\\
*={}  &   *={}   &   *={} \ar @{-} [l]  &   *={} \ar @{-} [l]  \\
*={}  &   *={} \ar @{-} [l] \ar@{-} [u] &   *={} \ar @{-} [l] \ar@{-} [u] &   *={} \ar @{-} [l] \ar@{-} [u] \\
*={} \ar @{-} [u] &   *={} \ar @{-} [l] \ar@{-} [u] &   *={} \ar @{-} [l] \ar@{-} [u] \\
*={} \ar @{-} [u] &   *={} \ar @{-} [l] \ar@{-} [u] }\endxy . $$
A \emph{column strict tableaux of shape} $\nu/\mu$ is a filling of $\nu/\mu$ satisfying (a) and (b) above.

\vspace{.5cm}

\noindent\textbf{Symmetric Functions.}  The symmetric group $S_n$ acts on the infinite set of variables $\{x_1,x_2,\ldots\}$ by permuting the indices $\leq n$ and fixing those $>n$.  Let
\begin{equation}\label{SymF}\Lambda_\C(x)=\{f\in \C[[x_1,x_2,\ldots]]\ :\ w(f)=f, \text{  permutations } w \}\nonumber\end{equation}
be \emph{the ring of symmetric functions in the variables $\{x_1,x_2,\ldots\}$}.  Let 
$$e_r(x)=\sum_{1 \leq i_1 < i_2 <\cdots<i_r} x_{i_1}x_{i_2}\cdots x_{i_r}\qquad \textnormal{and}\qquad
p_s(x)=\sum_{1\leq i} x_i^s, \qquad r,s\in \Z_{\geq 0},$$
be the $r$th \emph{elementary symmetric function} and the $s$th \emph{power sum symmetric function}, respectively.  For a partition $\nu=(\nu_1,\nu_2,\ldots, \nu_\ell)\in\mcP$, let
$$e_\nu(x)=e_{\nu_1}(x)e_{\nu_2}(x)\cdots e_{\nu_\ell}(x),\qquad p_\nu(x)=p_{\nu_1}(x)p_{\nu_2}(x)\cdots p_{\nu_\ell}(x)$$
and let
\begin{equation}\label{SF} s_\nu(x)=\det(e_{\nu_i^\prime-i+j}(x))\end{equation}
be the \emph{Schur function corresponding to} $\nu$.  The ring
\begin{equation}\label{SFbasis} \Lambda_\C(x)=\C\spanning\{e_\nu(x)\}=\C\spanning\{p_\nu(x)\}=\C\spanning\{s_\nu(x)\}, \end{equation}
and Pieri's rule says that if $\nu\in \mcP$, then
\begin{equation}\label{pieri} s_\nu(x) s_{(n)}(x)=\sum_{ \sh(P) = \gamma/\nu\atop \wt(P)=(n)} s_\gamma(x). \qquad \cite[\text{I}.5.16]{Mac}
\end{equation}
 
For each $t\in \C$ and partition $\nu$, let $P_\nu(x;t)$ denote the \emph{Hall-Littlewood symmetric function} \cite[\text{III}.2]{Mac}.  Since a precise definition is not necessary for this paper, it suffices to remark that
$$P_\nu(x;0)=s_\nu(x), \qquad  P_{(1^n)}(x;t)=e_n(x)$$
and for each $t\in \C$
$$\Lambda_\C(x)=\C\spanning\{P_\nu(x;t)\}.$$
(For additional details, see \cite[Chapter I]{Mac} on symmetric functions and \cite[Chapter III]{Mac} on Hall-Littlewood functions). 

\vspace{.15cm}

\noindent\textbf{Remark.}  It is usually sufficient to let $\{x_1,x_2,\ldots, x_K\}$ be a finite variable set (for $K$ much bigger than $|G_n|$); think of symmetric functions as polynomials rather than formal power series by setting $x_j=0$ for $j>K$ in the definitions above.  While Theorem \ref{ch} requires the infinite definition, I urge the reader to think in terms of the finite version everywhere else.

\vspace{.5cm}

\noindent\textbf{RSK correspondence.}   The classical RSK correspondence provides a combinatorial proof of the identity
\begin{equation}\label{2line} \prod_{i,j>0} \frac{1}{1-x_iy_j}=\sum_{\nu\vdash n, n\geq 0} s_\nu(x)s_\nu(y) \qquad \cite{Kn}\notag\end{equation}
by constructing a bijection between the matrices $b\in M_\ell(\Z_{\geq 0})$ and the set of pairs $(P(b),Q(b))$ of column strict tableaux with the same shape.  The bijection is as follows. 

If $P$ is a column strict tableau and $j\in \Z_{>0}$, let $P\leftarrow j$ be the column strict tableau given by the following algorithm
\begin{enumerate}
\item[(a)] Insert $j$ into the the first column of $P_k$ by displacing the smallest number $\geq j$.  If all numbers are $<j$, then place $j$ at the bottom of the first column.
\item[(b)] Iterate this insertion by inserting the displaced entry into the next column.
\item[(c)] Stop when the insertion does not displace an entry. 
\end{enumerate}

A \emph{two-line array} $\left(\begin{array}{cccc} i_1 & i_2 & \cdots & i_n\\ j_1 & j_2 & \cdots & j_n\end{array}\right)$ is a two-rowed array with $i_1\leq i_2\leq \cdots \leq i_n$ and $j_k\geq j_{k+1}$ if $i_k=i_{k+1}$.  If $b\in M_\ell(\Z_{\geq 0})$, then let $\vec{b}$ be the two-line array with $b_{ij}$ pairs $\tbinom{i}{j}$. 

For $b\in M_\ell(\Z_{\geq 0})$, suppose 
$$ \vec{b}=\left(\begin{array}{cccc} i_1 & i_2 & \cdots & i_n\\ j_1 & j_2 & \cdots & j_n\end{array}\right).$$
Then the pair $(P(b),Q(b))$ is the final pair in the sequence 
$$(\emptyset,\emptyset)=(P_0,Q_0),(P_1,Q_1),(P_2,Q_2),\ldots , (P_n,Q_n)=(P(b),Q(b)),$$
where $(P_k,Q_k)$ is a pair of column strict tableaux with the same shape given by
$$P_k=P_{k-1}\leftarrow j_k\qquad \textnormal{ and }\qquad \begin{array}{p{6.5cm}}  $Q_k$ is defined by $\sh(Q_k)=\sh(P_k)$ with $i_k$ in the new box $\sh(Q_k)/\sh(Q_{k-1})$.\end{array}$$
For example,
$$b=\left(\begin{array}{ccc} 1 & 1 & 0\\ 0 & 0 & 2\\ 0 & 1 & 0\end{array} \right)\qquad \text{corresponds to} \qquad\vec{b}=\left(\begin{array}{ccccc} 1 & 1 & 2 & 2 & 3 \\ 2 & 1 & 3 & 3 & 2\end{array} \right)$$
and provides the sequence
$$(\emptyset,\emptyset),
(\xy<0cm,.2cm>\xymatrix@R=.4cm@C=.4cm{ *={} & *={} \ar @{-} [l] \\ *={}\ar @{-} [u] & *={} \ar @{-} [l] \ar @{-} [u] \ar @{} [ul]|{2}}\endxy
,\xy<0cm,.2cm>\xymatrix@R=.4cm@C=.4cm{ *={} & *={} \ar @{-} [l] \\ *={}\ar @{-} [u] & *={} \ar @{-} [l] \ar @{-} [u] \ar @{} [ul]|{1}}\endxy)
, 
(\xy<0cm,.2cm>\xymatrix@R=.4cm@C=.4cm{ *={} & *={} \ar @{-} [l] & *={} \ar @{-} [l]\\ *={}\ar @{-} [u] & *={} \ar @{-} [l] \ar @{-} [u] \ar @{} [ul]|{1} & *={} \ar @{-} [l] \ar @{-} [u] \ar @{} [ul]|{2}}\endxy 
,\xy<0cm,.2cm>\xymatrix@R=.4cm@C=.4cm{ *={} & *={} \ar @{-} [l] & *={} \ar @{-} [l]\\ *={}\ar @{-} [u] & *={} \ar @{-} [l] \ar @{-} [u] \ar @{} [ul]|{1} & *={} \ar @{-} [l] \ar @{-} [u] \ar @{} [ul]|{1}}\endxy )
, 
\left(\xy<0cm,.4cm>\xymatrix@R=.4cm@C=.4cm{ *={} & *={} \ar @{-} [l] & *={} \ar @{-} [l]\\ *={}\ar @{-} [u] & *={} \ar @{-} [l] \ar @{-} [u] \ar @{} [ul]|{1} & *={} \ar @{-} [l] \ar @{-} [u] \ar @{} [ul]|{2}\\ *={}\ar @{-} [u] & *={} \ar @{-} [l] \ar @{-} [u] \ar @{} [ul]|{3}}\endxy ,
\xy<0cm,.4cm>\xymatrix@R=.4cm@C=.4cm{ *={} & *={} \ar @{-} [l] & *={} \ar @{-} [l]\\ *={}\ar @{-} [u] & *={} \ar @{-} [l] \ar @{-} [u] \ar @{} [ul]|{1} & *={} \ar @{-} [l] \ar @{-} [u] \ar @{} [ul]|{1}\\
*={}\ar @{-} [u] & *={} \ar @{-} [l] \ar @{-} [u] \ar @{} [ul]|{2}}\endxy\right)
 , 
\left(\xy<0cm,.4cm>\xymatrix@R=.4cm@C=.4cm{ *={} & *={} \ar @{-} [l] & *={} \ar @{-} [l]\\ *={}\ar @{-} [u] & *={} \ar @{-} [l] \ar @{-} [u] \ar @{} [ul]|{1} & *={} \ar @{-} [l] \ar @{-} [u] \ar @{} [ul]|{2}\\ *={}\ar @{-} [u] & *={} \ar @{-} [l] \ar @{-} [u] \ar @{} [ul]|{3} & *={} \ar @{-} [l] \ar @{-} [u] \ar @{} [ul]|{3}}\endxy
 ,\xy<0cm,.4cm>\xymatrix@R=.4cm@C=.4cm{ *={} & *={} \ar @{-} [l] & *={} \ar @{-} [l]\\ *={}\ar @{-} [u] & *={} \ar @{-} [l] \ar @{-} [u] \ar @{} [ul]|{1} & *={} \ar @{-} [l] \ar @{-} [u] \ar @{} [ul]|{1}\\ *={}\ar @{-} [u] & *={} \ar @{-} [l] \ar @{-} [u] \ar @{} [ul]|{2} & *={} \ar @{-} [l] \ar @{-} [u] \ar @{} [ul]|{2}}\endxy \right)
,
\left(\xy<0cm,.4cm>\xymatrix@R=.4cm@C=.4cm{ *={} & *={} \ar @{-} [l] & *={} \ar @{-} [l] & *={} \ar @{-} [l]\\ *={}\ar @{-} [u] & *={} \ar @{-} [l] \ar @{-} [u] \ar @{} [ul]|{1} & *={} \ar @{-} [l] \ar @{-} [u] \ar @{} [ul]|{2} & *={} \ar @{-} [l] \ar @{-} [u] \ar @{} [ul]|{3}\\ *={}\ar @{-} [u] & *={} \ar @{-} [l] \ar @{-} [u] \ar @{} [ul]|{2} & *={} \ar @{-} [l] \ar @{-} [u] \ar @{} [ul]|{3}}\endxy,
\xy<0cm,.4cm>\xymatrix@R=.4cm@C=.4cm{ *={} & *={} \ar @{-} [l] & *={} \ar @{-} [l] & *={} \ar @{-} [l]\\ *={}\ar @{-} [u] & *={} \ar @{-} [l] \ar @{-} [u] \ar @{} [ul]|{1} & *={} \ar @{-} [l] \ar @{-} [u] \ar @{} [ul]|{1} & *={} \ar @{-} [l] \ar @{-} [u] \ar @{} [ul]|{3}\\ *={}\ar @{-} [u] & *={} \ar @{-} [l] \ar @{-} [u] \ar @{} [ul]|{2} & *={} \ar @{-} [l] \ar @{-} [u] \ar @{} [ul]|{2}}\endxy\right)
$$
so that 
$$(P(b),Q(b))=\left(\xy<0cm,.4cm>\xymatrix@R=.4cm@C=.4cm{ *={} & *={} \ar @{-} [l] & *={} \ar @{-} [l] & *={} \ar @{-} [l]\\ *={}\ar @{-} [u] & *={} \ar @{-} [l] \ar @{-} [u] \ar @{} [ul]|{1} & *={} \ar @{-} [l] \ar @{-} [u] \ar @{} [ul]|{2} & *={} \ar @{-} [l] \ar @{-} [u] \ar @{} [ul]|{3}\\ *={}\ar @{-} [u] & *={} \ar @{-} [l] \ar @{-} [u] \ar @{} [ul]|{2} & *={} \ar @{-} [l] \ar @{-} [u] \ar @{} [ul]|{3}}\endxy,
\xy<0cm,.4cm>\xymatrix@R=.4cm@C=.4cm{ *={} & *={} \ar @{-} [l] & *={} \ar @{-} [l] & *={} \ar @{-} [l]\\ *={}\ar @{-} [u] & *={} \ar @{-} [l] \ar @{-} [u] \ar @{} [ul]|{1} & *={} \ar @{-} [l] \ar @{-} [u] \ar @{} [ul]|{1} & *={} \ar @{-} [l] \ar @{-} [u] \ar @{} [ul]|{3}\\ *={}\ar @{-} [u] & *={} \ar @{-} [l] \ar @{-} [u] \ar @{} [ul]|{2} & *={} \ar @{-} [l] \ar @{-} [u] \ar @{} [ul]|{2}}\endxy\right).$$

\vspace{.5cm}

\noindent \textbf{The general linear group.}  Let $G=\GL_n(\F_q)$, where $\F_q$ is the finite field with $q$ elements.  Let
$$ U=\left\{\left(\begin{array}{cccc} 1& * & \cdots & *\\ 0 & 1 & \ddots & \vdots \\ \vdots & \ddots & \ddots & *\\ 0 & \cdots & 0 & 1\end{array}\right)\right\}\subseteq G$$
be the subgroup of  unipotent, upper-triangular matrices.  For $1\leq i\neq j\leq n$, let $x_{ij}(t)$ be the matrix with ones on the diagonal, $t$ in the $(i,j)$th position and zeroes elsewhere.  Then
$$U=\langle x_{ij}(t)\ :\ 1\leq i<j\leq n\rangle.$$

The group $G$ has a double-coset decomposition given by
\begin{equation} \label{N} G=\bigsqcup_{v\in N} UvU,\qquad \text{where}\quad 
N=\left\{\begin{array}{p{5.2cm}} $n\times n$ matrices with entries from $\F_q$ and exactly one nonzero entry in each row and column\end{array}\right\}.\end{equation}
If $T\subseteq N$ is the subgroup of diagonal matrices and $W\subseteq N$ is the the subgroup of permutation matrices, then $N=WT$ and $TU=\mathbf{N}_G(U)$ is the normalizer of $U$ in $G$.  If necessary, specify the size of the group by a subscript such as $G_n$, $U_n$, $W_n$, etc.  Let $w_{(k)}\in W_k$ be the $k\times k$ matrix 
\begin{equation}\label{w_k} (w_{(k)})_{ij}=\delta_{j,n-i+1}. \qquad \textnormal{For example, } w_{(3)}=\left(\begin{array}{ccc} 0 & 0 & 1\\ 0 & 1 & 0\\  1 & 0 & 0\end{array}\right). \end{equation}

For $\mu=(\mu_1,\mu_2,\ldots,\mu_\ell)\models n$,   
\begin{equation}\label{Pmu} P_\mu=\left\{\left(\begin{array}{c@{}c@{}c@{}c} \begin{array}{|c|} \hline g_1\\ \hline\end{array}   & & & *\\
 & \begin{array}{|c|} \hline g_2\\ \hline\end{array} & & \\ & & \ddots & \\
 0 & & & \begin{array}{|c|} \hline g_\ell\\ \hline\end{array}\end{array}\right)\ :\ g_i\in G_{\mu_i}=\GL_{\mu_i}(\F_q)\right\}\end{equation}
is a \emph{parabolic subgroup} of $G$.  The \emph{Levi subgroup} and the \emph{unipotent radical} of $P_\mu$ are
\begin{equation}\label{LUmu} L_\mu=\left\{\left(\begin{array}{c@{}c@{}c@{}c} \begin{array}{|c|} \hline g_1\\ \hline\end{array}   & & & 0\\
 & \begin{array}{|c|} \hline g_2\\ \hline\end{array} & & \\ & & \ddots & \\
 0 & & & \begin{array}{|c|} \hline g_\ell\\ \hline\end{array}\end{array}\right)\ :\ g_i\in G_{\mu_i}\right\}
\ \text{and} \ 
U_\mu=\left\{\left(\begin{array}{c@{}c@{}c@{}c} \begin{array}{|c|} \hline Id_{\mu_1}\\ \hline\end{array}   & & & *\\
 & \begin{array}{|c|} \hline Id_{\mu_2}\\ \hline\end{array} & & \\ & & \ddots & \\
 0 & & & \begin{array}{|c|} \hline Id_{\mu_\ell}\\ \hline\end{array}\end{array}\right)\right\},\end{equation}
respectively, where $Id_k$ is the $k\times k$ identity matrix.  Note that $P_\mu=L_\mu U_\mu$ and $P_\mu=\mathbf{N}_G(U_\mu)$.

\end{section}

\begin{section}{An indexing for the standard basis of $\H_\mu$}

Let $G=\GL_n(\F_q)$.  Fix a nontrivial character $\psi : \F_q^+\rightarrow \C^*$ of the additive group of $\F_q$.  Let $\mu\models n$ and $\B_\mu$ be as in (\ref{Bmu}).  Since $x_{ij}(t)\in [U,U]$ for all $j>i+1$, the map $\psi_\mu: U\rightarrow \C^*$, defined by
\begin{equation}\label{psimu}\psi_\mu(x_{ij}(t))=\left\{ \begin{array}{ll} \psi(t), & \textnormal{if } j=i+1,i\notin \B_\mu,\\ 1, & \textnormal{otherwise,}\end{array}\right.\end{equation}
is a linear character of $U$.  Let
\begin{equation} \label{Hmu} \H_\mu=\End_G\left(\Ind_U^G(\psi_\mu)\right)\cong e_\mu \C G e_\mu,\qquad \text{where}\quad e_\mu=\frac{1}{|U|}\sum_{u\in U} \psi_\mu(u^{-1}) u.\end{equation}

The classical examples of unipotent Hecke algebras are the Yokonuma algebra $\H_{(1^n)}$ \cite{Y2} and the Gelfand-Graev Hecke algebra $\H_{(n)}$ \cite{St}.  A fundamental result is
\begin{theorem}[\cite{GG},\cite{Y1},\cite{St}]\label{HnC} For all $ n > 0$, $\H_{(n)}$ is commutative. 
\end{theorem} 
\noindent This theorem will follow from Theorem \ref{Hn}.

An analysis of (\ref{Tbasis}) implies that $\{T_v= e_\mu ve_\mu\ :\ v\in N_\mu\}$ is a basis for $\H_\mu$ , where
\begin{equation} N_\mu=\{v\in N\ :\ u,vuv^{-1}\in U \textnormal{ implies } \psi_\mu(u)=\psi_\mu(vuv^{-1})\}.\qquad \cite[(11.30)]{CR}\end{equation}

Suppose $a\in M_\ell(\F_q[X])$ is an $\ell\times \ell$ matrix with polynomial entries.  Let $d(a_{ij})$ be the degree of the polynomial $a_{ij}$.  Define the \emph{degree row sums} and the \emph{degree column sums} of $a$ to be the compositions 
$$d^\rightarrow(a)=(d^\rightarrow(a)_1,d^\rightarrow(a)_2,\ldots, d^\rightarrow(a)_\ell)\qquad \textnormal{and}\qquad d^\downarrow(a)=(d^\downarrow(a)_1,d^\downarrow(a)_2,\ldots, d^\downarrow(a)_\ell),$$
where
$$d^\rightarrow(a)_i=\sum_{j=1}^\ell d(a_{ij})\qquad \textnormal{and} \qquad d^\downarrow(a)_j=\sum_{i=1}^\ell d(a_{ij}).$$
Let 
\begin{equation} M_\mu=\{a\in M_{\ell(\mu)}(\F_q[X])\ :\ d^\rightarrow(a)=d^\downarrow(a)=\mu, a_{ij}\textnormal{ monic, } a_{ij}(0)\neq 0\}.\end{equation}
For example, 
$$\begin{array}{c@{} c}\left(\begin{array}{cccc} X+1 & 1 & 1 & X+2\\ X+3 & X^3+2X+3 & 1 & X+2\\ 1 & X^2+4X+2 & X+2 & 1\\
1 & 1 & X^2+3X+1 & X^2+2\end{array}\right) & \begin{array}{c} \ss{(1+0+0+1=2)}\\ \ss{(1+3+0+1=5)}\\ \ss{(0+2+1+0=3)}\\ \ss{(0+0+2+2=4)}\vspace{.05cm}\end{array}\\
\begin{array}{c@{\hspace{.3cm}}c@{\hspace{.4cm}}c@{\hspace{.3cm}}c} \ss{(1+1+0+0=2)} & \ss{(0+3+2+0=5)} & \ss{(0+0+1+2=3)} & \ss{(1+1+0+2=4)}\end{array} & \end{array} \in M_{(2,5,3,4)}.$$

Suppose $(f)=(a_0+a_1X^{i_1}+a_2X^{i_2}+\cdots+a_rX^{i_r}+X^n)\in M_{(n)}$ is a $1\times 1$ matrix with $a_0,a_1,\ldots,a_r\neq 0$.  Let
\begin{equation}\label{ftov_f} v_{(f)}=w_{(n)}(a_0w_{(i_1)}\oplus a_1 w_{(i_2-i_1)}\oplus\cdots \oplus a_r w_{(n-i_r)})\in N,\end{equation}
where $w_{(k)}\in W_k$ is as in (\ref{w_k}) and by convention $v_{(1)}$ is the matrix with height and width 0.  For example,
$$v_{(a +bX^3+ cX^4 +X^6)}=
\left(\begin{array}{cccccc} 
0 & 0 & 0 & 0 & 0 & 1 \\
0 & 0 & 0 & 0 & 1 & 0 \\
0 & 0 & 0 & 1 & 0 & 0 \\
0 & 0 & 1 & 0 & 0 & 0 \\ 
0 & 1 & 0 & 0 & 0 & 0 \\ 
1 & 0 & 0 & 0 & 0 & 0 \end{array}\right)
\left(\begin{array}{c @{}c@{}c}\begin{array}{|ccc|}\hline 0 & 0 & a  \\
0 & a & 0  \\
a & 0 & 0 \\ \hline \end{array} & \begin{array}{c} 0 \\ 0\\ 0\end{array} & \begin{array}{cc} 0 & 0 \\ 0 & 0 \\ 0 & 0\end{array}\\
\begin{array}{ccc} 0 & 0 & 0\end{array} & \begin{array}{|c|}\hline b\\ \hline\end{array} & \begin{array}{cc} 0 & 0\end{array}\\
\begin{array}{ccc} 0 & 0 & 0 \\ 0 & 0 & 0\end{array} & \begin{array}{c} 0\\ 0\end{array} & \begin{array}{|cc|} \hline 0 & c \\ c & 0 \\ \hline\end{array}  
\end{array}\right)=
\left(\begin{array}{cccccc} 
0 & 0 & 0 & 0 & c & 0 \\
0 & 0 & 0 & 0 & 0 & c \\
0 & 0 & 0 & b & 0 & 0 \\ 
a & 0 & 0 & 0 & 0 & 0 \\ 
0 & a & 0 & 0 & 0 & 0 \\ 
0 & 0 & a & 0 & 0 & 0 \end{array}\right).$$

For each $a\in M_\mu$ construct a matrix $v_a\in N$ by partitioning the rows and columns of an $n\times n$ matrix according to $\mu=(\mu_1,\mu_2,\ldots,\mu_\ell)\models n$ and setting
\begin{equation} \label{v_a} v_a=\begin{array}{@{}c@{}l@{}} \begin{array}{c@{\hspace{.1cm}}c@{\hspace{.1cm}}c@{\hspace{.1cm}}c}\overbrace{\xymatrix@R=0cm@C=3.3cm{*={}\ar@{} [r] & *={}}}^{\mu_1\ \mathrm{columns}} & \overbrace{\xymatrix@R=0cm@C=3.3cm{*={}\ar@{} [r] & *={}}}^{\mu_2\ \mathrm{columns}} & \overbrace{\xymatrix@R=0cm@C=3.3cm{*={}\ar@{} [r] & *={}}}^{\mu_3\ \mathrm{columns}} & \xymatrix@R=0cm@C=1.6cm{*={}\ar @{} [r]|{\dd{\cdots}} & *={}}\end{array} &  \vspace{-.1cm}\\
\left(\xy<0cm,4.3cm>\xymatrix@R=.6cm@C=.85cm{ 
*={} & *={} & *={} & *={} & *={} \ar @{} [ddddllll]|{0} \ar @{-} [dddddddddddddd] & *={} & *={} & *={} & *={} \ar @{} [dddllll]|{0} \ar @{-} [dddddddddddddd] & *={} & *={} & *={} & *={} \ar @{} [ddllll]|{0} \ar @{-} [dddddddddddddd] & *={}  & *={} \ar @{} [dll]|{\dd{\cdots}}\\
*={} & *={} & *={} & *={} & *={} & *={} & *={} & *={} & *={} & *={} & *={} & *={} & *={} \ar @{} [dl]|{v_{(a_{13})}} \ar @{-} [l] \ar @{--} [rr] & *={} & *={} \ar @{} [dll]|{0} \\
*={} & *={} & *={} & *={} & *={} & *={} & *={} & *={} & *={} \ar @{} [dl]|{v_{(a_{12})}} \ar @{-} [l] \ar @{--} [rrrrrr]  & *={} & *={} & *={} \ar@{-} [u] \ar @{-} [r] \ar @{--} [dddddddddddd] \ar @{} [dlll]|{0} & *={} \ar @{} [dl]|{0} & *={}  & *={} \ar @{} [dll]|{0}\\
*={} & *={} & *={} & *={} & *={} \ar @{} [dl]|{v_{(a_{11})}} \ar @{-} [l] \ar @{--} [rrrrrrrrrr] & *={} & *={} & *={} \ar@{-} [u] \ar @{-} [r] \ar @{--} [ddddddddddd] \ar @{} [dlll]|{0} & *={}\ar @{} [dl]|{0} & *={} & *={} & *={}\ar @{} [dlll]|{0} & *={} \ar @{} [dl]|{0} & *={} & *={} \ar @{} [dll]|{0}\\
*={} \ar @{-} [rrrrrrrrrrrrrr] & *={} & *={} & *={} \ar@{-} [u] \ar @{-} [r] \ar @{--} [dddddddddd] \ar @{} [ddddlll]|{0} & *={}  \ar @{} [dddl]|{0} & *={} & *={} & *={} \ar @{} [dddlll]|{0} & *={}\ar @{} [ddl]|{0} & *={} & *={} & *={} \ar @{} [ddlll]|{0} & *={} \ar @{} [dl]|{0}& *={} & *={} \ar @{} [dll]|{\dd{\cdots}} \\
*={} & *={} & *={} & *={} & *={} & *={} & *={} & *={} & *={} & *={} & *={} & *={} \ar @{} [dl]|{v_{(a_{23})}} \ar @{-} [l] \ar @{-} [d] \ar @{--} [rrr]  & *={} \ar @{} [dl]|{0} & *={}  & *={} \ar @{} [dll]|{0} \\
*={} & *={} & *={} & *={} & *={} & *={} & *={} & *={}  \ar @{} [dl]|{v_{(a_{22})}} \ar @{-} [l] \ar @{-} [d] \ar @{--} [rrrrrrr]   & *={} \ar @{} [dl]|{0} & *={} & *={} \ar@{-} [u] \ar @{-} [r] \ar @{--} [dddddddd] \ar @{} [dll]|{0}  & *={} \ar @{} [dl]|{0} & *={} \ar @{} [dl]|{0} & *={}  & *={} \ar @{} [dll]|{0} \\
*={} & *={} & *={} & *={}  \ar @{} [dl]|{v_{(a_{21})}} \ar @{-} [l] \ar @{-} [d] \ar @{--} [rrrrrrrrrrr]   & *={} \ar @{} [dl]|{0} & *={} & *={} \ar@{-} [u] \ar @{-} [r] \ar @{--} [ddddddd] \ar @{} [dll]|{0} & *={} \ar @{} [dl]|{0} & *={} \ar @{} [dl]|{0} & *={} & *={} \ar @{} [dll]|{0} & *={} \ar @{} [dl]|{0} & *={} \ar @{} [dl]|{0} & *={} & *={} \ar @{} [dll]|{0}\\
*={} \ar @{-} [rrrrrrrrrrrrrr] & *={} & *={}  \ar@{-} [u] \ar @{-} [r] \ar @{--} [dddddd] \ar @{} [ddddll]|{0} & *={} \ar @{} [dddl]|{0} & *={} \ar @{} [dddl]|{0} & *={} & *={} \ar @{} [ddll]|{0} & *={} \ar @{} [ddl]|{0} & *={} \ar @{} [ddl]|{0} & *={}\ar @{} [dl]|{0} & *={}  & *={} \ar @{} [dl]|{0} & *={} \ar @{} [dl]|{0} & *={} & *={} \ar @{} [dll]|{\dd{\cdots}}  \\
*={} & *={} & *={} & *={} & *={} & *={} & *={} & *={} & *={} & *={} & *={}  \ar @{} [dl]|{v_{(a_{33})}} \ar @{-} [l] \ar @{-} [d] \ar @{--} [rrrr] & *={} \ar @{} [dl]|{0} & *={} \ar @{} [dl]|{0} & *={} & *={} \ar @{} [dll]|{0}\\
*={} & *={} & *={} & *={} & *={} & *={} & *={} \ar @{} [dl]|{v_{(a_{32})}} \ar @{-} [l] \ar @{-} [d] \ar @{--} [rrrrrrrr] & *={} \ar @{} [dl]|{0} & *={} \ar @{} [dl]|{0} & *={} \ar@{-} [u] \ar @{-} [r] \ar @{--} [dddd] \ar @{} [dl]|{0}  & *={}\ar @{} [dl]|{0} & *={}\ar @{} [dl]|{0} & *={}\ar @{} [dl]|{0} & *={} & *={} \ar @{} [dll]|{0} \\
*={} & *={} & *={} \ar @{} [dl]|{v_{(a_{31})}} \ar @{-} [l] \ar @{-} [d] \ar @{--} [rrrrrrrrrrrr]   & *={}\ar @{} [dl]|{0} & *={}\ar @{} [dl]|{0} & *={} \ar@{-} [u] \ar @{-} [r] \ar @{--} [ddd]\ar @{} [dl]|{0} & *={}\ar @{} [dl]|{0} & *={}\ar @{} [dl]|{0} & *={}\ar @{} [dl]|{0} & *={}\ar @{} [dl]|{0} & *={} \ar @{} [dl]|{0}& *={}\ar @{} [dl]|{0} & *={}\ar @{} [dl]|{0} & *={} & *={}\ar @{} [dll]|{0} \\
*={} \ar @{-} [rrrrrrrrrrrrrr] & *={} \ar@{-} [u] \ar @{-} [r] \ar @{--} [dd]  \ar @{} [ddl]|{\dd{\vdots}} & *={}\ar @{} [ddl]|{0} & *={}\ar @{} [ddl]|{0} & *={} \ar @{} [ddl]|{0} & *={}  \ar @{} [ddl]|{\dd{\vdots}} & *={} \ar @{} [ddl]|{0}& *={}\ar @{} [ddl]|{0} & *={}\ar @{} [ddl]|{0} & *={} \ar @{} [ddl]|{\dd{\vdots}} & *={}\ar @{} [ddl]|{0} & *={} \ar @{} [ddl]|{0}& *={} \ar @{} [ddl]|{0}& *={} & *={} \ar @{} [ddll]|{\dd{\ddots}} \\
*={} & *={} & *={} & *={} & *={} & *={} & *={} & *={} & *={} & *={} & *={} & *={} & *={} & *={} & *={} \\
*={} & *={} & *={} & *={} & *={} & *={} & *={} & *={} & *={} & *={} & *={} & *={} & *={} & *={} & *={}}
\endxy \right) & \begin{array}{@{}l@{}}\left. \xy<0cm,1.2cm> \xymatrix@R=2.4cm@C=0cm{*={}\ar @{} [d] \\ *={}}\endxy\hspace{-.1cm} \right\} \hspace{-.05cm} \rotatebox[origin=1cm]{-90}{$\ss{\mu_1\ \mathrm{rows}}$}\\ \left.\xy<0cm,1.2cm> \xymatrix@R=2.4cm@C=0cm{*={}\ar @{} [d] \\ *={}}\endxy\hspace{-.1cm} \right\} \hspace{-.05cm} \rotatebox[origin=1cm]{-90}{$\ss{\mu_2\ \mathrm{rows}}$} \\ \left.\xy<0cm,1.2cm> \xymatrix@R=2.4cm@C=0cm{*={}\ar @{} [d] \\ *={}}\endxy\hspace{-.1cm} \right\} \hspace{-.05cm} \rotatebox[origin=1cm]{-90}{$\ss{\mu_3\ \mathrm{rows}}$}\\ \xy<0cm,.7cm>\xymatrix@R=1.35cm@C=0cm{*={} \ar @{} [d]|{\dd{\vdots}}\\ *={}}\endxy \end{array} \end{array}\ , \end{equation}
where $a_{ij}$ is the $(i,j)$th entry of $a$ and $v_{(a_{ij})}$ is as in (\ref{ftov_f}). 

\begin{theorem} \label{basis} The map
$$\begin{array}{rcl} M_\mu & \longrightarrow & N_\mu\\
a & \mapsto & v_a,\end{array}$$
given by (\ref{v_a}) is a bijection.
\end{theorem}

\noindent\textbf{Remarks.}  When $\mu=(n)$ this theorem says that the map $(f)\mapsto v_{(f)}$ of (\ref{ftov_f}) is a bijection between $M_{(n)}$ and $N_{(n)}$.  

\begin{proof}
Using the remark following the theorem, it is straightforward to reconstruct $a$ from $v_a$. Therefore the map is invertible, and it suffices to show 
\begin{enumerate}
\item[(a)]  the map is well-defined ($v_a\in N_\mu$),
\item[(b)]  the map is surjective.
\end{enumerate}
To show (a) and (b), we investigate the matrices $N_\mu$.  Suppose $v\in N_\mu$.  Let
$$\begin{array}{rcl} v_i & = & \textnormal{ the nonzero entry in the } i\textnormal{th column of }  v, \vspace{.15cm}\\
v(i) & = & \textnormal{the row number of the nonzero entry in the } i\textnormal{th column of } v,\end{array}$$
so that $\{v_1,v_2,\ldots,v_n\}$ are the nonzero entries of $v$ and $(v(1),v(2),\ldots, v(n))$ is the permutation determined by setting all the nonzero entries of $v$ to 1.  By (\ref{psimu}),
\begin{equation} \psi_\mu(x_{ij}(t))=\left\{\begin{array}{ll} \psi(t), & \textnormal{if } j=i+1 \textnormal{ and } i\notin \B_\mu,\\ 1, & \textnormal{otherwise.}\end{array}\right. \tag{A} \end{equation} 
Recall that $v\in N_\mu$ if and only if $u,vuv^{-1}\in U$ implies $\psi_\mu(u)=\psi_\mu(vuv^{-1})$. That is, $v\in N_\mu$ if and only if for all $1\leq i<j\leq n$ such that $v(i)<v(j)$,
\begin{align}  
\psi_\mu(x_{ij}(t)) & =  \psi_\mu(vx_{ij}(t)v^{-1}) \notag\\
& =   \psi_\mu(x_{v(i)v(j)}(v_i t v_j^{-1})) \notag\\
& = \left\{\begin{array}{ll} \psi(v_it v_j^{-1}), & \textnormal{if } v(j)=v(i)+1 \textnormal{ and } v(i)\notin \B_\mu,\\ 1, & \textnormal{otherwise.}\end{array}\right. \tag{B}
 \end{align}
Compare (A) and (B) to obtain that $v\in N_\mu$ if and only if for all $1\leq i<j\leq n$ such that $v(i)<v(j)$,
\begin{enumerate}
\item[(i)] If $i\notin \B_\mu$ and $v(i)\in \B_\mu$, then $j\neq i+1$,
\item[(ii)] If $i\in \B_\mu$ and $v(i)\notin \B_\mu$, then $v(j)\neq v(i)+1$,
\item[(iii)] If $i,v(i)\notin \B_\mu$, then $j=i+1$ if and only if $v(j)=v(i)+1$,
\item[$\mathrm{(iii)}^\prime$]  If $i,v(i)\notin \B_\mu$ and $v(j)=v(i)+1$, then $v_i=v_{i+1}$.
\end{enumerate}
We can visualize the implications of the conditions (i)$\mathrm{-(iii)}^\prime$ in the following way.  Partition the rows and columns of $v\in N_\mu$ by $\mu$.  For example, $\mu=(2,3,1)$ partitions $v$ according to
$$\left(\begin{array}{cc|ccc|c} * & * & * & * & * & *\\
* & * & * & * & * & *\\ \hline
* & * & * & * & * & *\\
* & * & * & * & * & *\\
* & * & * & * & * & *\\ \hline
* & * & * & * & * & *\end{array} \right).$$
Suppose the nonzero entry of $v$ in column $i$ is above a horizontal line but not next to a vertical line.  Then condition (i) implies that $v(i+1)<v(i)$, so
\begin{equation} \begin{array}{c @{\hspace{1cm}} c} \begin{array}{ccccc} 
 & & 0 & *  &\\
 & & 0 & *  & \\
 & & 0 & * & \\
0 & 0 & a & 0  & 0\\ \hline\end{array} & \begin{array}{p{6cm}} where  $a$ is the nonzero entry and  $*$  indicates possible locations for nonzero entries in the next column.\end{array}\end{array} \tag{I}\end{equation}
Similarly, condition (ii) implies
\begin{equation}\begin{array}{cccc|}  & & & 0\\
 & & &  0\\
 0 & 0  & 0 & a\\
 * & * & * & 0\\
 & & & 0\end{array} \tag{II}\end{equation}
and conditions (iii) and $\mathrm{(iii)}^\prime$ imply
\begin{equation}\begin{array}{ccccc}
& & & 0 & *  \\ 
& & & 0 & *  \\
& & & 0 & *   \\
0 & 0 & 0 & a & 0 \\
* & * & * & 0 & 0   \end{array}\qquad \textnormal{or}\qquad 
\begin{array}{cccccc} 
  & &  & 0 & 0 & \\
 & &  & 0 & 0 & \\
 0 & 0 & 0  & a & 0 & 0  \\
0 &  0 & 0 & 0 & a &  0 \\
 &      &   & 0 & 0 &
\end{array}. \tag{III}\end{equation}
In the case $\mu=(n)$ condition (III) implies that every $v\in N_{(n)}$ is of the form 
$$\left( \begin{array}{c@{}c@{}c@{}c} 0 &  & \hspace{.5cm} & \begin{array}{|c|} \hline a_r Id_{i_r}\\ \hline\end{array} \\  &  & \xy<0cm,.3cm>\xymatrix@R=.08cm@C=.09cm{*={} & *={} & *={\cdot} \\ *={} & *={\cdot} & *={}\\ *={\cdot} & *={}  & *={}}\endxy &  \\
 & \begin{array}{|c|} \hline  a_2 Id_{i_2} \\ \hline\end{array} &  & \\
\begin{array}{|c|} \hline a_1 Id_{i_1}\\ \hline\end{array} &  &  & 0\end{array}\right)=w_{(n)}(a_1w_{(i_1)}\oplus a_2 w_{(i_2)}\oplus \cdots \oplus a_r w_{(i_r)}),$$
where $(i_1,i_2,\ldots,i_r)\models n$, $a_i\in \F_q^*$, and $w_{(k)}\in W_k$ is as in (\ref{w_k}).  In fact, this observation proves that the map $(f)\mapsto v_{(f)}$ is a bijection between $M_{(n)}$ and $N_{(n)}$ (mentioned in the remark).

Note that since the matrices $v_a$ satisfy (I), (II) and (III), $v_a\in N_\mu$, proving (a).  On the other hand, (I), (II) and (III) imply that each $v\in N_\mu$  must be of the form $v=v_a$ for some $a\in M_\mu$, proving (b).  
\end{proof}
\end{section}

\begin{section}{An RSK-insertion via the representation theory of $\H_\mu$}

Let $\mathcal{S}$ be a set.  An \emph{$\mathcal{S}$-partition} $\lambda=(\lambda^{(s_1)}, \lambda^{(s_2)},\ldots )$ is a sequence of partitions indexed by the elements of $\mathcal{S}$.    Let
\begin{equation} \label{P^S} \mcP^\mathcal{S}=\{ \mathcal{S}\text{-partitions}\}.\end{equation}
The following discussion defines two sets $\Theta$ and $\Phi$, so that $\Theta$-partitions index the irreducible characters of $G$ and $\Phi$-partitions index the conjugacy classes of $G$.

Let $L_n=\Hom(\F_{q^n}^*,\C^*)$ be the character group of $\F_{q^n}^*$.  If $\gamma\in L_m$, then let
$$\begin{array}{rcl} \dd{\gamma_{(r)}: \F_{q^{mr}}^*} &  \longrightarrow  & \dd{\C^*}\vspace{.15cm}\\ 
x & \mapsto & \gamma(x^{1+q^r+q^{2r}+\cdots+q^{m(r-1)}})\end{array}$$
Thus if $n=mr$, then we may view $L_m\subseteq L_n$ by identifying $\gamma\in L_m$ with $\gamma_{(r)}\in L_n$.  Define
$$L=\bigcup_{n\geq 0} L_n.$$
The \emph{Frobenius maps} are
$$\begin{array}{rcl} F:\bar\F_q & \rightarrow & \bar\F_q\\ x\ & \mapsto & x^q \end{array} \qquad\textnormal{and}\qquad \begin{array}{rcl} F:L & \rightarrow & L\\ \gamma\ & \mapsto & \gamma^q \end{array},$$
where $\bar\F_q$ is the algebraic closure of $\F_q$.  

The map
$$\begin{array}{rcl} \{$F$\text{-orbits of } \bar\F_q^*\} & \longrightarrow & \{f\in \F_q[t]\ :\ f\text{ is monic, irreducible, and } f(0)\neq 0\}\\
\dd{\{x,x^q,x^{q^2},\ldots,x^{q^{k-1}}\}} & \mapsto & \dd{f_x=\prod_{i=1}^{k-1} (t-x^{q^i}), \qquad \text{where } x^{q^k}=x\in \bar\F_q^*}\end{array}$$
is a bijection such that the size of the $F$-orbit of $x$ equals the degree $d(f_x)$ of $f_x$.   Let 
\begin{equation}\label{PhiTheta} \Phi=\left\{f\in \C[t]\ :\ \begin{array}{p{3.6cm}} $f$ is monic, irreducible and  $f(0)\neq 0$\end{array}\right\}\quad \text{and}\quad  \Theta=\{F\textnormal{-orbits in } L\}.\end{equation} 
If $\eta$ is a $\Phi$-partition and $\lambda$ is a $\Theta$-partition, then let
$$|\eta |=\sum_{f\in \Phi}d(f)|\eta^{(f)}| \qquad \textnormal{and}\qquad |\lambda |=\sum_{\vphi\in \Theta}|\vphi||\lambda^{(\vphi)}|$$
be the \emph{size} of $\eta$ and $\lambda$, respectively.  Let the sets $\mcP^\Phi$ and $\mcP^\Theta$ be as in (\ref{P^S}) and let
\begin{equation}\label{P_n} \mcP_n^\Phi=\{\eta\in \mcP^\Phi\ :\ |\eta|=n\}\qquad\text{and}\qquad \mcP_n^\Theta=\{ \lambda\in \mcP^\Theta \ :\ |\lambda|=n\}.\end{equation}
\begin{theorem}[Green \cite{Gr}] \label{Green} Let $G_n=\GL_n(\F_q)$.
\begin{enumerate}
\item[(a)] $\mcP_n^\Phi$ indexes the conjugacy classes $K^\eta$ of $G_n$,
\item[(b)] $\mcP_n^\Theta$ indexes the irreducible $G_n$-modules $G_n^\lambda$.
\end{enumerate}
\end{theorem}
Suppose $\lambda\in \mcP^\Theta$.  A \emph{column strict tableau} $P=(P^{(\vphi_1)},P^{(\vphi_2)},\ldots)$ \emph{of shape $\lambda$} is a column strict filling of $\lambda$ by positive integers.  That is, $P^{(\vphi)}$ is a column strict tableau of shape $\lambda^{(\vphi)}$.  Write $\sh(P)=\lambda$.  The \emph{weight of $P$} is the composition $\wt(P)=(\wt(P)_1,\wt(P)_2,\ldots )$ given by
$$\wt(P)_i=\sum_{\vphi\in \Theta} |\vphi| \left( \begin{array}{c} \textnormal{number of}\\ i \textnormal{ in } P^{(\vphi)}\end{array}\right).$$

If $\lambda\in \mcP^\Theta$ and $\mu$ is a composition, then let
\begin{equation}\label{Hmuhat^l} \hat{\H}_\mu^\lambda=\{\textnormal{column strict tableaux } P\ :\ \sh(P)=\lambda, \wt(P)=\mu\}\end{equation}
and 
\begin{equation} \label{Hmuhat}\hat\H_\mu=\{\lambda\in \mcP^\Theta\ :\ \hat\H_\mu^\lambda \textnormal{ is not empty}\}.\end{equation}

The following theorem is a consequence of (\ref{SW}) and a theorem proved by Zelevinsky \cite{Z} (see Theorem \ref{ZelG}).  A proof of Zelevinsky's theorem is in Section \ref{Zproof}.

\begin{theorem}\label{Zel}  The set $\hat{\H}_\mu$ indexes the irreducible $\H_\mu$-modules $\H_\mu^\lambda$ and
$$\dim(\H_\mu^\lambda)= |\hat\H_\mu^\lambda |.$$
\end{theorem}

The $(\H_\mu,\H_\mu)$-bimodule decomposition 
$$\H_\mu\cong \bigoplus_{\lambda\in \hat\H_\mu} \H_\mu^\lambda \otimes \H_\mu^\lambda\quad \text{implies}\quad |N_\mu|=\dim(\H_\mu)=\bigoplus_{\lambda\in \hat\H_\mu} \dim(\H_\mu^\lambda)^2 =\sum_{\lambda\in \hat\H_\mu} |\hat\H_\mu^\lambda|^2.$$
Theorem \ref{RSK}, below, gives a combinatorial proof of this identity.

Encode each matrix $a\in M_\mu$ as a $\Phi$-sequence
$$(a^{(f_1)}, a^{(f_2)},\ldots),\qquad f_i\in \Phi,$$
where $a^{(f)}\in M_{\ell(\mu)}(\Z_{\geq 0})$ is given by
$$a_{ij}^{(f)}= \textnormal{highest power of } f \textnormal{ dividing } a_{ij}.$$ 
Note that this is an entry by entry ``factorization" of $a$ such that
$$a_{ij}=\prod_{f\in \Phi} f^{a_{ij}^{(f)}}.$$
Recall from page \pageref{2line} the classical RSK correspondence
$$\begin{array}{rcl} M_\ell(\Z_{\geq 0}) & \longrightarrow & \left\{\begin{array}{p{5cm}} Pairs $(P,Q)$ of column strict tableaux of the same shape\end{array}\right\}\vspace{.15cm}\\
b &  \mapsto   & (P(b),Q(b)).\end{array}$$

\begin{theorem} \label{RSK}  For $a\in M_\mu$, let $P(a)$ and $Q(a)$ be the $\Phi$-column strict tableaux given by 
$$P(a)=(P(a^{(f_1)}), P(a^{(f_2)}),\ldots)\qquad \textnormal{and}\qquad Q(a)=(Q(a^{(f_1)}), Q(a^{(f_2)}),\ldots) \qquad \textnormal{for } f_i\in \Phi.$$
Then the map
$$\begin{array}{rcccl} N_\mu &\longrightarrow & M_\mu & \longrightarrow & \left\{\begin{array}{p{4.5cm}} Pairs $(P,Q)$ of $\Phi$-column strict tableaux of the same shape and weight $\mu$\end{array}\right\}\vspace{.15cm}\\
v & \mapsto & a_v & \mapsto & (P(a_v),Q(a_v)),\end{array}$$
is a bijection, where the first map is the inverse of the bijection of Theorem \ref{basis}.
\end{theorem}
By the construction above, the map is well-defined and since all the steps are invertible, the map is a bijection. 

For example, suppose $\mu=(7,5,3,2)$ and $f,g,h\in \Phi$ are such that $d(f)=1$, $d(g)=2$, and $d(h)=3$. Then
$$a_v=\left(\begin{array}{cccc} g &  f^2 h & 1 & 1 \\ h & 1 & g & 1\\ 1 & 1 & f & f^2\\ g & 1 & 1 & 1 \end{array}\right)\in 
M_{\xymatrix@R=.1cm@C=.1cm{*={} & *={} \ar @{-} [l] &  *={} \ar @{-} [l] & *={} \ar @{-} [l] & *={} \ar @{-} [l] & *={} \ar @{-} [l] & *={} \ar @{-} [l] & *={} \ar @{-} [l] \\
*={} \ar @{-} [u] & *={} \ar @{-} [l] \ar @{-} [u] &  *={} \ar @{-} [l] \ar @{-} [u] & *={} \ar @{-} [l] \ar @{-} [u] & *={} \ar @{-} [l] \ar @{-} [u] & *={} \ar @{-} [l] \ar @{-} [u] & *={} \ar @{-} [l] \ar @{-} [u] & *={} \ar @{-} [l] \ar @{-} [u] \\
*={} \ar @{-} [u] & *={} \ar @{-} [l] \ar @{-} [u] &  *={} \ar @{-} [l] \ar @{-} [u] & *={} \ar @{-} [l] \ar @{-} [u] & *={} \ar @{-} [l] \ar @{-} [u] & *={} \ar @{-} [l] \ar @{-} [u]\\
*={} \ar @{-} [u] & *={} \ar @{-} [l] \ar @{-} [u] &  *={} \ar @{-} [l] \ar @{-} [u] & *={} \ar @{-} [l] \ar @{-} [u]\\
*={} \ar @{-} [u] & *={} \ar @{-} [l] \ar @{-} [u] &  *={} \ar @{-} [l] \ar @{-} [u]}}
$$
corresponds to the sequence
$$(a_v^{(f_1)},a_v^{(f_2)},\ldots)=\left(\left(\begin{array}{cccc}0 &  2 & 0 & 0\\ 0 & 0 & 0 & 0\\ 0 & 0 & 1 & 2\\ 0 & 0 & 0 & 0 \end{array}\right)^{(f)}, 
\left(\begin{array}{cccc}1 &  0 & 0 & 0\\ 0 & 0 & 1 & 0\\ 0 & 0 & 0 & 0\\ 1 & 0 & 0 & 0 \end{array}\right)^{(g)}, 
\left(\begin{array}{cccc}0 &  1 & 0 & 0\\ 1 & 0 & 0 & 0\\ 0 & 0 & 0 & 0\\ 0 & 0 & 0 & 0\end{array}\right)^{(h)}\right)$$
and
$$(P(a_v),Q(a_v))=\left(\xy<0cm,.4cm>\xymatrix@R=.4cm@C=.4cm{
	*={} & *={} \ar @{-} [l] & *={} \ar @{-} [l] & *={} \ar @{-} [l]  \\
	*={} \ar @{-} [u] &   *={} \ar @{-} [l] \ar@{-} [u] \ar@{} [ul] |{2} &   *={} \ar @{-} [l] \ar@{-} [u] \ar@{} [ul] |{2} & *={} \ar @{-} [l] \ar@{-} [u] \ar@{} [ul] |{4} \\
*={} \ar @{-} [u] &   *={} \ar @{-} [l] \ar@{-} [u] \ar@{} [ul] |{3} & *={} \ar @{-} [l] \ar@{-} [u] \ar@{} [ul] |{4} }\endxy^{(f)},\ 
\xy<0cm,.4cm>\xymatrix@R=.4cm@C=.4cm{
	*={} & *={} \ar @{-} [l] & *={} \ar @{-} [l]  \\
	*={} \ar @{-} [u] &   *={} \ar @{-} [l] \ar@{-} [u] \ar@{} [ul] |{1}&   *={} \ar @{-} [l] \ar@{-} [u] \ar@{} [ul] |{1} \\ *={} \ar @{-} [u] &   *={} \ar @{-} [l] \ar@{-} [u] \ar@{} [ul] |{3}}\endxy^{(g)}, \  \xy<0cm,.4cm>\xymatrix@R=.4cm@C=.4cm{
	*={} & *={} \ar @{-} [l] & *={} \ar @{-} [l] \\
	*={} \ar @{-} [u] &   *={} \ar @{-} [l] \ar@{-} [u] \ar@{} [ul] |{1} &  *={} \ar @{-} [l] \ar@{-} [u] \ar@{} [ul] |{2}}\endxy^{(h)}\right)
\left(\xy<0cm,.4cm>\xymatrix@R=.4cm@C=.4cm{
	*={} & *={} \ar @{-} [l] & *={} \ar @{-} [l] & *={} \ar @{-} [l]  \\
	*={} \ar @{-} [u] &   *={} \ar @{-} [l] \ar@{-} [u] \ar@{} [ul] |{1}&   *={} \ar @{-} [l] \ar@{-} [u] \ar@{} [ul] |{1} & *={} \ar @{-} [l] \ar@{-} [u] \ar@{} [ul] |{3} \\
*={} \ar @{-} [u] &   *={} \ar @{-} [l] \ar@{-} [u] \ar@{} [ul] |{3} & *={} \ar @{-} [l] \ar@{-} [u] \ar@{} [ul] |{3}}\endxy^{(f)},\ 
\xy<0cm,.4cm>\xymatrix@R=.4cm@C=.4cm{
	*={} & *={} \ar @{-} [l] & *={} \ar @{-} [l] & *={}  \\
	*={} \ar @{-} [u] &   *={} \ar @{-} [l] \ar@{-} [u] \ar@{} [ul] |{1}&   *={} \ar @{-} [l] \ar@{-} [u] \ar@{} [ul] |{4} \\ *={} \ar @{-} [u] &  *={} \ar @{-} [l] \ar@{-} [u] \ar@{} [ul] |{2}}\endxy^{(g)}, \  \xy<0cm,.4cm>\xymatrix@R=.4cm@C=.4cm{
	*={} & *={} \ar @{-} [l] & *={} \ar @{-} [l] \\
	*={} \ar @{-} [u] &   *={} \ar @{-} [l] \ar@{-} [u] \ar@{} [ul] |{1} &  *={} \ar @{-} [l] \ar@{-} [u] \ar@{} [ul] |{2}}\endxy^{(h)}\right).$$

\end{section}

\begin{section}{Zelevinsky's decomposition of $\Ind_U^G(\psi_\mu)$}\label{Zproof}  

This section proves the theorem 
\begin{theorem}[Zelevinsky \cite{Z}] \label{ZelG}  Let $U$ be the subgroup of unipotent upper-triangular matrices of $G=\GL_n(\F_q)$, $\mu\models n$ and $\psi_\mu$ be as in (\ref{psimu}).
$$\Ind_U^G(\psi_\mu)=\bigoplus_{\lambda\in \hat\H_\mu} (G^\lambda)^{\oplus |\hat\H_\mu^\lambda|}.$$
\end{theorem}
Theorem \ref{Zel} follows from this theorem and double-centralizer theory (\ref{SW}).  The following will
\begin{enumerate}
\item[(i)] establish the necessary connection between symmetric functions and the representation theory of $G$,
\item[(ii)] prove Theorem \ref{ZelG} for the case when $\ell(\mu)=1$,
\item[(iii)] generalize to arbitrary $\mu$.
\end{enumerate}
The proof below uses the ideas of Zelevinsky's proof, but explicitly uses symmetric functions to prove the results.  Specifically,  the following discussion through the proof of Theorem \ref{Hn} corresponds to \cite[Sections 9-11]{Z} and Theorem \ref{ZelG} corresponds to \cite[Theorem 12.1]{Z}.

\vspace{.5cm}

\noindent\textbf{Preliminaries to the proof.}  Suppose $\lambda\in \mcP^\Theta$ and $\eta\in \mcP^\Phi$ (see (\ref{P_n}).   Let $\chi^\lambda$ be the irreducible character corresponding to the irreducible $G$-module $G^\lambda$ and let $\kappa^\eta$ be the characteristic function corresponding to the conjugacy class $K^\eta$ (see Theorem \ref{Green}), given by
$$\kappa^\eta(g)=\left\{\begin{array}{ll} 1, & \textnormal{if } g\in K^\eta,\\ 0, & \textnormal{otherwise,} \end{array}\right. \qquad \text{for } g\in G_{|\eta|}.$$ 
Define
$$R=\C\spanning\{\chi^\lambda\ :\ \lambda\in \mcP^\Theta\}=\C\spanning\{\kappa^\eta\ :\  \eta\in \mcP^\Phi\}.$$ 
The space $R$ has an inner product defined by
$$\langle \chi^\lambda,\chi^\nu\rangle=\delta_{\lambda\nu},$$
and multiplication
\begin{equation}\label{Rmult} \chi^\lambda\circ \chi^\nu=\Indf_{L_{(r,s)}}^{G_{r+s}}(\chi^\lambda\otimes \chi^\nu)=\Ind_{P_{(r,s)}}^{G_{r+s}}\left(\Inf_{L_{(r,s)}}^{P_{(r,s)}}(\chi^\lambda\otimes \chi^\nu)\right),\qquad \textnormal{for } \lambda\in\mcP^\Theta_r, \nu\in \mcP^\Theta_s,\end{equation}
where if $p\in P_\mu=L_\mu U_\mu$ decomposes as $p=lu$ for $u\in U_\mu$, $l\in L_\mu$, then
$$\Inf_{L_\mu}^{P_\mu}(\chi^{\gamma_1}\otimes  \cdots \otimes \chi^{\gamma_\ell})(p)=\chi^{\gamma_1}(l_1) \cdots  \chi^{\gamma_\ell}(l_\ell), \qquad \textnormal{for } l=(l_1\oplus\cdots\oplus l_\ell), l_i\in G_{\mu_i}, \gamma_i\in \mcP^\Theta_{\mu_i}.$$

For each $\vphi\in \Theta$, let $\{Y_1^{(\vphi)}, Y^{(\vphi)}_2,\ldots\}$ be an infinite set of variables, and let
$$\Lambda_\C=\bigotimes_{\vphi\in \Theta}\Lambda_\C(Y^{(\vphi)}),$$
where $\Lambda_\C(Y^{(\vphi)})$ is the ring of symmetric functions in $\{Y^{(\vphi)}_1,Y^{(\vphi)}_2,\ldots\}$ (see page \pageref{SymF}).  For each $f\in \Phi$, define an additional set of variables $\{X_1^{(f)}, X_2^{(f)},\ldots \}$ such that the symmetric functions in the $Y$ variables are related to the symmetric functions in the $X$ variables by the transform
\begin{equation}\label{transform}p_k(Y^{(\vphi)})=(-1)^{k|\vphi|-1}\sum_{x\in {\F_{q^{k|\vphi|}}}^*} \xi(x) p_{\frac{k|\vphi|}{d(f_x)}}(X^{(f_x)}),\end{equation}
where $\xi\in \vphi$, $f_x\in \Phi$ is the irreducible polynomial that has $x$ as a root, and $p_{\frac{a}{b}}(X^{(f)})=0$ if $\frac{a}{b}\notin\Z_{\geq 0}$.  Then
$$\Lambda_\C=\bigotimes_{f\in \Phi} \Lambda_\C(X^{(f)}).$$

For $\nu\in \mcP$, let $s_\nu(Y^{(\vphi)})$ be the Schur function and $P_\nu(X^{(f)};t)$ be the Hall-Littlewood symmetric function (see page \pageref{SymF}).  Define
\begin{equation}\label{sP} s_\lambda=\prod_{\vphi\in \Theta} s_{\lambda^{(\vphi)}}(Y^{(\vphi)})\qquad \textnormal{and}\qquad P_\eta=q^{-n(\eta)}\prod_{f\in \Phi} P_{\eta^{(f)}}(X^{(f)}; q^{-d(f)}),\end{equation}
where $n(\eta)=\sum_{f\in \Phi} d(f) n(\eta^{(f)})$ and $n(\mu)=\sum_{i=1}^\ell (i-1)\mu_i$, for $\mu$ a composition.
The ring 
$$\Lambda_\C=\C\spanning\{s_\lambda\ :\ \lambda\in \mcP^\Theta\}=\C\spanning\{P_\eta\ :\ \eta\in \mcP^\Phi\}$$
has an inner product given by
$$\langle s_\lambda,s_\nu\rangle=\delta_{\lambda\nu}.$$
\begin{theorem}[Green \cite{Gr},Macdonald \cite{Mac}] \label{ch} The linear map
$$\begin{array}{rcl}\ch: R & \longrightarrow & \Lambda_\C\\
\chi^\lambda & \mapsto & s_\lambda,\qquad \text{for } \lambda\in \mcP^\Theta \\
\kappa^\eta & \mapsto & P_\eta,\qquad \text{for } \eta\in \mcP^\Phi, \end{array}$$
is an algebra isomorphism that preserves the inner product.
\end{theorem}

The unipotent conjugacy classes $K^\eta$ satisfy $\eta^{(f)}=\emptyset$ unless $f=t-1$.  Let
$$\U=\C\spanning\{\kappa^\eta\ :\ \eta^{(f)}=\emptyset,\textnormal{ unless } f=t-1\}\subseteq R$$
be the subalgebra of unipotent class functions.  Note that by (\ref{sP}) and Theorem \ref{ch}
$$\ch(\U)=\Lambda_\C(X^{(t-1)}).$$
Consider the projection $\pi: R\rightarrow \U$ which is an algebra homomorphism given by
$$(\pi\chi^\lambda)(g)=\left\{\begin{array}{ll} \chi^\lambda(g), & \textnormal{if } g\in G \textnormal{ is unipotent,}\\ 0, & \textnormal{otherwise}, \end{array}\right. \qquad \lambda\in \mcP^\Theta.$$
Then $\tilde\pi=\pi\circ\ch^{-1}:\Lambda_\C\rightarrow \C$ is given by
\begin{equation}\label{tpi} \begin{split}\tilde\pi(p_k(Y^{(\vphi)})) & =  \tilde\pi\left((-1)^{k|\vphi|-1}\sum_{x\in {\F_{q^{k|\vphi|}}}^*} \xi(x) p_{\frac{k|\vphi|}{d(f_x)}}(X^{(f_x)})\right) \quad \textnormal{(by (\ref{transform}))} \\
& =  (-1)^{k|\vphi|-1}\xi(1)p_{\frac{k|\vphi|}{1}}(X^{(t-1)})+0\\
& =  (-1)^{k|\vphi|-1} p_{k|\vphi|}(X^{(t-1)}).\end{split}\end{equation}

\vspace{.5cm}

\noindent\textbf{The decomposition of $\Ind_U^G(\psi_{(n)})$.}  The representation $\Ind_U^G(\psi_{(n)})$ is the Gelfand-Graev module, and with (\ref{SW}), Theorem \ref{Hn} proves that $\H_{(n)}$ is commutative.

For $\lambda\in \mcP^\Theta$, let
$$\hgt(\lambda)=\max\{\ell(\lambda^{(\vphi)})\ :\ \vphi\in \Theta\}.$$

\begin{theorem}\label{Hn}
$$\ch(\Ind_U^G(\psi_{(n)}))=\sum_{\lambda\in \mcP^\Theta_n,\hgt(\lambda)=1} s_\lambda.$$
\end{theorem}

\begin{proof}  Let 
\begin{equation}\label{Psi} \begin{array}{rcl} \Psi: R & \longrightarrow & \C\\
\chi^\lambda & \mapsto & \langle \chi^\lambda, \Ind_U^G(\psi_{(n)})\rangle\end{array} \qquad \text{and} \qquad \tilde\Psi: \Lambda_\C\overset{\ch^{-1}}{\longrightarrow} R \overset{\Psi}{\longrightarrow} \C.\end{equation}  For any group $H$, let $1_H$ be the trivial character of $H$, $e_H=\frac{1}{|H|}\sum_{h\in H}h$, and $\langle \chi, \gamma\rangle_H=\frac{1}{|H|}\sum_{h\in H} \chi(h)\gamma(h^{-1})$, for all class functions $\gamma$ and $\chi$ of $H$.

The proof is in six steps.
\begin{enumerate}
\item[(a)] $\tilde\Psi(e_k(Y^{(1)}))=\delta_{k1}$, where $1$ is the trivial character of $\F_q^*$,
\item[(b)] $\Psi(\chi^\lambda)=\dim(e_{(n)}G^\lambda)$ for $\lambda\in \mcP^\Theta$,
\item[(c)] $\tilde\Psi(fg)=\tilde\Psi(f)\tilde\Psi(g)$ for all $f,g\in \Lambda_\C(Y^{(1)})$, where $1$ is the trivial character of $\F_q^*$,
\item[(d)] $\Psi\circ \pi=\Psi$,
\item[(e)] $\tilde\Psi(f(Y^{(\vphi)}))=\tilde\Psi(f(Y^{(1)}))$ for all $f\in \Lambda_\C(Y^{(\vphi)})$,
\item[(f)] $\tilde\Psi(s_\lambda)=\delta_{\hgt(\lambda)1}$.
\end{enumerate}

(a) An argument similar to the argument in \cite[pgs. 285-286]{Mac} shows that $$\ch^{-1}(e_k(Y^{(1)}))=1_{G_k}$$
(see \cite[Theorem 4.9 (a)]{HR} for details).  Therefore, by Frobenius reciprocity and the orthogonality of characters, 
$$\tilde\Psi(e_k(Y^{(1)}))=\langle 1_{G_k},\Ind_U^G(\psi_{(n)})\rangle=\langle 1_{U_k}, \psi_{(n)}\rangle_{U_k} =\delta_{k1}.$$ 

(b) Since there exists an idempotent $e$ such that $G^\lambda\cong \C G e$ and $\Ind_U^G(\psi_{(n)})\cong \C G e_{(n)}$, the map
$$\begin{array}{rcl} e_{(n)} \C G e & \longrightarrow & \Hom_G(G^\lambda, \C G e_{(n)})\vspace{.15cm}\\
                     e_{(n)} g e & \mapsto & \begin{array}{rcl}\gamma_g: \C G e & \rightarrow & \C G e_{(n)}\\ xe & \mapsto & x e g e_{(n)}\end{array}\end{array}$$
is a vector space isomorphism (using an argument similar to the proof of \cite[(3.18)]{CR}).  Thus, $$\Psi(\chi^\lambda)=\langle \chi^\lambda,\Ind_U^G(\psi_{(n)})\rangle=\dim(\Hom_G(G^\lambda, \Ind_U^G(\psi_{(n)})))=\dim(e_{(n)} \C G e)=\dim(e_{(n)} G^\lambda).$$ 

(c)   By (a), $\tilde\Psi(e_r(Y^{(1)})) \tilde\Psi (e_s(Y^{(1)}))=\delta_{r1}\delta_{s1}$.  Since $\Lambda_\C(Y^{(1)})=\C[e_1(Y^{(1)}), e_2(Y^{(1)}),\ldots]$, it therefore suffices to show that
$$\tilde\Psi(e_r(Y^{(1)})e_s(Y^{(1)}))=\delta_{r1}\delta_{s1}.$$
Suppose $r+s=n$ and let $P=P_{(r,s)}$. Then
$$\tilde\Psi(e_r(Y^{(1)})e_s(Y^{(1)}))=\Psi(\Ind_{P}^{G_n}(1_P))=\dim(e_{(n)}\C G e_P).$$
Since $T\subseteq P$, $e_P=e_{(1^n)} e_P$, $G=\bigsqcup_{v\in N} UvU$, and $N=WT$,
$$e_{(n)}\C G e_P=e_{(n)}\C G e_{(1^n)} e_P=\C\spanning\{ e_{(n)}w e_{(1^n)}e_P\ :\ w\in W\}.$$
If there exists $1\leq i\leq n$ such that $w(i)=w(i)+1$, then
$$e_{(n)}w e_{(1^n)}=e_{(n)}w x_{i,i+1}(t)e_{(1^n)}=e_{(n)}x_{w(i), w(i)+1}(t) w e_{(1^n)}= \psi(t) e_{(n)}w e_{(1^n)}.$$
Therefore, $e_{(n)} w e_{(1^n)}=0$ unless $w=w_{(n)}$.  If $r>1$ of $s>1$, then there exists $1\leq i\leq n$ such that $x_{i+1,i}(t)\in P_{(r,s)}$, so
$$e_{(n)}w_{(n)} e_P=e_{(n)}w_{(n)} x_{i+1,i}(t) e_P=e_{(n)}x_{n-i,n-i+1}(t) w_{(n)} e_P=\psi(t) e_{(n)}w_{(n)} e_P=0.$$
In particular,
$$\dim(e_{(n)}\C G e_P)=0.$$
If $r=s=1$, then $P_{(1,1)}$ is upper-triangular, so
$$e_{(2)} w_{(2)} e_P\neq 0$$
and $\dim(e_{(2)} \C G e_P)=1$, giving $\tilde\Psi(e_r(Y^{(1)})e_s(Y^{(1)}))=\delta_{r1}\delta_{s1}$.

(d) By Frobenius reciprocity,
$$\langle \chi^\lambda, \Ind_{U_n}^{G_n}(\psi_{(n)})\rangle =  \langle \Res_{U_n}^{G_n}(\chi^\lambda), \psi_{(n)}\rangle_{U_n} =  \langle \Res_{U_n}^{G_n}(\pi(\chi^\lambda)), \psi_{(n)}\rangle_{U_n}= \langle \pi(\chi^\lambda), \Ind_{U_n}^{G_n}(\psi_{(n)})\rangle,$$
so $\Psi=\Psi\circ \pi$.

(e) Induct on $n$, using (c) and the identity
$$(-1)^{n-1}p_n(Y^{(1)})=ne_n(Y^{(1)})-\sum_{r=1}^{n-1} (-1)^{r-1}p_r(Y^{(1)}) e_{n-r}(Y^{(1)}),\qquad \cite[\text{I}.2.11^\prime]{Mac}$$
to obtain $\tilde\Psi(p_n(Y^{(1)}))=1$.  Note that
$$\tilde\Psi(p_n(Y^{(\vphi)}))=\tilde\Psi(\pi(p_n(Y^{(\vphi)})))= \tilde\Psi((-1)^{|\vphi|n-1}p_{|\vphi|n}(X^{(t-1)}))
=\tilde\Psi(\pi(p_{|\vphi|k}(Y^{(1)})))$$
$$=\tilde\Psi(p_{|\vphi| k}(Y^{(1)}))=1=\tilde\Psi(p_k(Y^{(1)})).$$
Since $\tilde\Psi$ is multiplicative  on $\Lambda_\C(Y^{(1)})$, 
$$\tilde\Psi(p_\nu(Y^{(\vphi)}))=1=\tilde\Psi(p_\nu(Y^{(1)})),\qquad \textnormal{for all partitions } \nu.$$
In particular, since $\tilde\Psi$ is linear and $\Lambda_\C(Y^{(\vphi)})=\C\spanning\{p_\nu(Y^{(\vphi)})\}$,
$$\tilde\Psi(f(Y^{(\vphi)}))=\tilde\Psi(f(Y^{(1)})),\qquad \textnormal{for all } f\in \Lambda_\C(Y^{(\vphi)}).$$
Note that (e) also implies that $\tilde\Psi$ is multiplicative on all of $\Lambda_\C$.

(f)  Note that
$$\tilde\Psi(s_\lambda)=\tilde\Psi\left(\prod_{\vphi\in \Theta} s_{\lambda(\vphi)}(Y^{(\vphi)})\right)=\tilde\Psi\left(\prod_{\vphi\in \Theta} s_{\lambda(\vphi)}(Y^{(1)})\right)=\prod_{\vphi\in \Theta} \tilde\Psi(s_{\lambda(\vphi)}(Y^{(1)})),
$$
where the last two equalities follow from (e) and (c), respectively. 
By definition $s_\nu(Y^{(1)})=\det(e_{\nu_i^\prime-i+j}(Y^{(1)}))$, so
$$\tilde\Psi(s_\nu(Y^{(1)}))=\left\{\begin{array}{ll} 1, & \textnormal{if } \ell(\nu)=1,\\ 0, & \textnormal{otherwise,}\end{array}\right.$$
implies
\begin{equation*} \ch(\Ind_U^G(\psi_{(n)}))=\sum_{\lambda\in \mcP^\Theta_n} \tilde\Psi(s_\lambda)s_\lambda=\sum_{\lambda\in \mcP^\Theta_n, \hgt(\lambda)=1} s_\lambda.\qedhere\end{equation*}
\end{proof}

\vspace{.5cm}

\noindent\textbf{Decomposition of $\Ind_U^G(\psi_\mu)$.}  Suppose $\lambda, \nu\in \mcP^\Theta$.  A \emph{column strict tableau $P$ of shape $\lambda$ and weight $\nu$} is a column strict filling of $\lambda$ such that for each $\vphi\in \Theta$,
$$\sh(P^{(\vphi)})=\lambda^{(\vphi)} \qquad \textnormal{ and }\qquad \wt(P^{(\vphi)})= \nu^{(\vphi)}.$$
We can now prove the theorem stated at the beginning of this section: 

\vspace{.25cm}

\noindent \textbf{Theorem \ref{ZelG} (\cite{Z})}  Let $U$ be the subgroup of unipotent upper-triangular matrices of $G=\GL_n(\F_q)$, $\mu\models n$ and $\psi_\mu$ be as in (\ref{psimu}).  Then
$$\Ind_U^G(\psi_\mu)=\bigoplus_{\lambda\in \hat\H_\mu} (G^\lambda)^{\oplus |\hat\H_\mu^\lambda|}.$$

\begin{proof}   Note that
$$\Ind_U^{P_\mu}(\psi_\mu)\cong \C P_\mu e_\mu = \C P_\mu e_{[\mu]} e_{[\mu]}^\prime,$$
where
\begin{equation} \label{e[mu]} e_{[\mu]}=\frac{1}{| U\cap L_\mu|}\sum_{u\in U\cap L_\mu} \psi_\mu(u^{-1}) u \qquad \textnormal{and} \qquad e_{[\mu]}^\prime=\frac{1}{|U_\mu|}\sum_{u\in U_\mu} u.\end{equation}
Thus,
$$\begin{array}{rcl} \dd{\Ind_U^{P_\mu}(\psi_\mu)} & \cong & \dd{\Inf_{L_\mu}^{P_\mu}\left(\Ind_{U\cap L_\mu}^{L_\mu}(\psi_\mu)\right)} \vspace{.15cm}\\
& \cong & \dd{\Inf_{L_\mu}^{P_\mu}\left(\Ind_{U_{\mu_1}}^{G_{\mu_1}}(\psi_{(\mu_1)})\otimes\Ind_{U_{\mu_2}}^{G_{\mu_2}}(\psi_{(\mu_2)})\otimes \cdots \otimes \Ind_{U_{\mu_\ell}}^{G_{\mu_\ell}}(\psi_{(\mu_\ell)}) \right)}\end{array}$$
In particular, by the definition of multiplication in $R$ (\ref{Rmult}),
$$\Gamma_\mu=\ch(\Ind_U^G(\psi_\mu))=\Gamma_{\mu_1}\Gamma_{\mu_2}\cdots \Gamma_{\mu_\ell},\qquad
\textnormal{where } \Gamma_{\mu_i}=\sum_{\lambda\in \mcP^\Theta_{\mu_i}, \hgt(\lambda)=1} s_\lambda.$$
Pieri's rule (\ref{pieri}) implies that for $\lambda\in \mcP^\Theta_r$, $\nu\in \mcP^\Theta_s$ and $\hgt(\nu)=1$, 
$$s_\lambda s_\nu=\sum_{\gamma\in \mcP^\Theta_{r+s}, |\hat\H_{\nu}^{\gamma/\lambda}|\neq 0} s_\gamma, \qquad \text{so}\qquad \Gamma_\mu=\sum_{\lambda\in \mcP^\Theta} K_{\lambda\mu} s_\lambda,$$
where 
$$\begin{array}{rcl} \dd{K_{\lambda\mu}} & = & \dd{\Card\{ \emptyset=\gamma_0\subset \gamma_1\subset \gamma_2\subset \cdots \subset \gamma_\ell=\lambda\ :\  |\hat\H_{(\mu_{i+1})}^{\gamma_{i+1}/\gamma_i}|=1\}}\vspace{.15cm}\\
& = & \Card\{ \textnormal{column strict tableaux of shape } \lambda \textnormal{ and weight } \mu\}
=|\hat\H_\mu^\lambda|.\end{array}$$
By Green's Theorem (Theorem \ref{ch}), $\ch$ is an isomorphism, so
\begin{equation*}\Ind_U^G(\psi_\mu)=\ch^{-1}(\Gamma_\mu)=\sum_{\lambda\in \hat\H_\mu}|\hat\H_\mu^\lambda| \ch^{-1}(s_\lambda)=\sum_{\lambda\in \hat\H_\mu} (G^\lambda)^{\oplus |\hat\H_\mu^\lambda|}.\qedhere\end{equation*} 
\end{proof}

\end{section}

\begin{section}{A weight space decomposition of $\H_\mu$-modules}

Let $\mu=(\mu_1,\mu_2,\ldots,\mu_\ell)\models n$ and let $P_\mu$, $L_\mu$ and $U_\mu$ be as in (\ref{Pmu}) and (\ref{LUmu}).  Recall that
$$e_\mu=\frac{1}{|U|}\sum_{u\in U}\psi_\mu(u^{-1}) u.$$

\begin{theorem}\label{cartan}   For $a\in M_\mu$, let $T_{a}=e_\mu v_a e_\mu$ with $v_a$ as in (\ref{v_a}).  Then the map
$$\begin{array}{rcl}\H_{(\mu_1)}\otimes \H_{(\mu_2)}\otimes \cdots\otimes \H_{(\mu_\ell)} & \longrightarrow & \H_\mu\\
T_{(f_1)}\otimes T_{(f_2)}\otimes \cdots \otimes T_{(f_\ell)} & \mapsto & T_{(f_1)\oplus (f_2)\oplus \cdots \oplus (f_\ell)}, \qquad \textnormal{for } (f_i)\in M_{(\mu_i)}\end{array}$$
is an injective algebra homomorphism with image $\L_\mu=e_\mu P_\mu e_\mu=e_\mu L_\mu e_\mu$.
\end{theorem}
\begin{proof}  Note that
$$T_{(f_1)}\otimes T_{(f_2)}\otimes \cdots \otimes T_{(f_\ell)}=\frac{1}{|U\cap L_\mu|^2}\sum_{x_i,y_i\in U_{\mu_i}} \left(\prod_{i=1}^\ell \psi_\mu(x_i^{-1}y_i^{-1})\right) x_1v_{(f_1)}y_1\otimes x_2v_{(f_2)} y_2\otimes \cdots \otimes x_\ell v_{(f_\ell)} y_\ell.$$
Since $U=(L_\mu\cap U)(U_\mu)$, $L_\mu\cap U\cong U_{\mu_1}\times U_{\mu_2}\times\cdots \times U_{\mu_\ell}$, and $\psi_\mu$ is trivial on $U_\mu$,
$$\begin{array}{rcl}\dd{T_{(f_1)\oplus (f_2)\oplus \cdots \oplus (f_\ell)}} & = & \dd{\frac{1}{|U|^2}\sum_{x,y\in U}\psi_\mu(x^{-1}y^{-1}) x(v_{(f_1)}\oplus v_{(f_2)}\oplus \cdots \oplus v_{(f_\ell)})y}\vspace{.15cm}\\
& = & 
\dd{\frac{1}{|U\cap L_\mu|^2}\sum_{x_i,y_i\in U_{\mu_i}}\psi_\mu(x_1^{-1}y_1^{-1}\oplus\cdots\oplus x_\ell^{-1}y_\ell^{-1})e_{[\mu]}^\prime x_1v_{(f_1)} y_1 \oplus \cdots \oplus x_\ell v_{(f_\ell)} y_\ell e_{[\mu]}^\prime},\end{array}$$
where $e_{[\mu]}^\prime$ is as in (\ref{e[mu]}).
Since $L_\mu\subseteq \mathbf{N}_G(U_\mu)$, the idempotent $e_{[\mu]}^\prime$ commutes with $x_1 v_{(f_1)} y_1 \oplus \cdots \oplus x_\ell v_{(f_\ell)} y_\ell$ and
$$T_{(f_1)\oplus (f_2)\oplus \cdots \oplus (f_\ell)}=\frac{e_{[\mu]}^\prime}{|L\cap U|^2}\sum_{x_i,y_i\in U_{\mu_i}}\left(\prod_{i=1}^\ell\psi_\mu(x_i^{-1}y_i^{-1})\right) x_1 v_{(f_1)} y_1\oplus \cdots \oplus x_\ell v_{(f_\ell)} y_\ell.$$
Consequently, the map multiplies by $e_{[\mu]}^\prime$ and changes $\otimes$ to $\oplus$, so it is an algebra homomorphism.  Since the map sends basis elements to basis elements, it is also injective. 
\end{proof}

Let $\L_\mu$ be as in Theorem \ref{cartan}. By Theorem \ref{HnC} each $\H_{(\mu_i)}$ is commutative,  so $\L_\mu$ is commutative and all the irreducible $\L_\mu$-modules $\L_\mu^\gamma$ are one-dimensional.  Theorem \ref{Zel} implies that
$$\hat\H_{(\mu_i)}=\{\Theta\text{-partitions } \lambda\ :\  |\lambda|=\mu_i, \hgt(\lambda)=1\}.$$ 
indexes the irreducible $\H_{(\mu_i)}$-modules.  Therefore, the set
\begin{equation}\label{Lmuhat}\hat\L_\mu=\hat\H_{(\mu_1)}\times\hat\H_{(\mu_2)}\times \cdots\times \hat\H_{(\mu_\ell)}=\{\gamma=(\gamma_1,\gamma_2,\ldots, \gamma_\ell)\ :\ \gamma_i\in \hat\H_{(\mu_i)}\} \end{equation}
indexes the irreducible $\L_\mu$-modules.  Identify $\gamma\in \hat\L_\mu$ with the map $\gamma:\L_\mu\rightarrow \C$ such that
$$yv=\gamma(y) v,\qquad \textnormal{ for all } y\in \L_\mu, v\in \L_\mu^\gamma.$$ 

For $\gamma\in \hat\L_\mu$, define the \emph{$\gamma$-weight space} $V_\gamma$ of an $\H_\mu$-module $V$ to be
$$V_\gamma=\{ v\in V\ :\ yv=\gamma(y) v, \textnormal{ for all } y\in \L_\mu\}.$$  
Then
$$V \cong \bigoplus_{\gamma\in \hat\L_\mu}V_\gamma.$$

Let $\lambda\in \mcP^\Theta$ and $\gamma\in \hat\L_\mu$.  A \emph{column strict tableau $P$ of shape $\lambda$ and weight $\gamma$} is column strict filling of $\lambda$ such that for each $\vphi\in \Theta$, 
$$\sh(P^{(\vphi)})=\lambda^{(\vphi)} \qquad \textnormal{and}\qquad \wt(P^{(\vphi)})=(|\gamma_1^{(\vphi)}|,|\gamma_2^{(\vphi)}|,\ldots, |\gamma_\ell^{(\vphi)}|),$$ 
where $|\gamma_i^{(\vphi)}|$ is the number of boxes in the partition $\gamma_i^{(\vphi)}$ (which has length 1).

\begin{theorem} Let $\H_\mu^\lambda$ be an irreducible $\H_\mu$-module and $\gamma\in \hat\L_\mu$.  Then
$$\dim(\H_\mu^\lambda)_\gamma=\Card\{\textnormal{column strict tableaux of shape } \lambda \textnormal{ and weight } \gamma\}=|\hat\H_\gamma^\lambda|.$$
\end{theorem}
\begin{proof}  By double-centralizer theory and Frobenius reciprocity,
$$\dim((\H_\mu^\lambda)_\gamma)=\langle\Res^{\H_\mu}_{\L_\mu}(\H_\mu^\lambda),\L_\mu^\gamma\rangle=\langle \Res^G_{P_\mu}(G^\lambda),P_\mu^\gamma\rangle=\langle G^\lambda, \Ind_{P_\mu}^G(P_\mu^\gamma)\rangle,$$
where $P_\mu^\gamma=\Inf_{L_\mu}^{P_\mu}(L_\mu^\gamma)$.  Therefore, 
$$\dim((\H_\mu^\lambda)_\gamma)=c_\gamma^\lambda,\qquad \textnormal{where}\quad s_{\gamma_1}s_{\gamma_2}\cdots s_{\gamma_\ell}=\sum_{\lambda\in \mcP^\Theta} c_\gamma^\lambda s_\lambda.$$
Pieri's rule (\ref{pieri}) implies $c_\gamma^\lambda=|\hat\H_\gamma^\lambda|$.
\end{proof}
\end{section}


\begin{thebibliography}{99}
\bibitem[Cu1]{CuS} Curtis, C. ``Representations of Hecke algebras." \emph{Ast\'erisque} \textbf{9} (1988): 13-60.
\bibitem[Cu2]{Cu} Curtis, C.  ``On the Gelfand-Graev Representations of a Reductive Group over a Finite Field."  \emph{Journal of Algebra}  \textbf{157} (1993): 517-533.
\bibitem[CR]{CR} Curtis, C. and Reiner, I.  \emph{Methods of Representation Theory, Vol. 1.}    New York: John Wiley and Sons,  1981.
\bibitem[CS]{CS} Curtis, C. and Shinoda, K.  ``Unitary Kloosterman Sums and the Gelfand-Graev Representation of $\GL_2^*$." \emph{Journal of Algebra} \textbf{216} (1999): 431-447.
\bibitem[Dr]{Dr} Drinfeld, V.  ``Quantum Groups." \emph{Proceedings of the International Congress of Mathematicians, Berkeley, California, 1-2 1986.}  Providence, RI: American Mathematical Society, 1987. 798-820.
\bibitem[GG]{GG} Gelfand, I.M. and Graev, M.I.  ``Construction of irreducible representations of simple algebraic groups over a finite field." \emph{Soviet Mathematics Doklady} \textbf{3} (1962): 1646-1649.
\bibitem[GW]{GW} Goodman R. and Wallach N.  \emph{Representations and invariants of the classical groups.} Encyclopedia of mathematics and its applications, 68.  Cambridge: Cambridge University Press, 1998.
\bibitem[Gr]{Gr} Green, J. A.  ``The Characters of the finite general linear groups."  \emph{Transactions of the American Mathematical Society} \textbf{80} (1955): 402-447.
\bibitem[HR]{HR} Halverson, T. and Ram, A. ``Bitraces for $\GL_n(\F_q)$ and the Iwahori-Hecke algebra of type $A_{n-1}$."   \emph{Indagationes Mathematicae} \textbf{10} (1999): 247-268.
\bibitem[Iw]{Iw} Iwahori, N.  ``On the structure of a Hecke ring of a Chevalley group over a finite field." \emph{Journal of the Faculty of Science, University of Tokyo}  \textbf{10} (1964): 215-236.
\bibitem[IM]{IM} Iwahori, N. and Matsumoto, H.  ``On some Bruhat decomposition and the structure of Hecke rings of $p$-adic Chevalley groups."  \emph{Institut des Hautes \'Etudes Scientifiques, Publications Math\'ematiques}  \textbf{25} (1965): 5-48.
\bibitem[Ji]{Ji} Jimbo, M.  ``A $q$-analogue of $U(\mathfrak{g}\mathfrak{l}(N+1))$, Hecke algebra, and the Yang-Baxter equation." \emph{Letters in  Mathematical Physics} \textbf{11} (1986): 247-252.
\bibitem[Jo1]{Jo1} Jones, V.  ``Index for subfactors." \emph{Inventiones Mathematicae} \textbf{72} (1983):  1-25.
\bibitem[Jo2]{Jo2} Jones, V.  ``Hecke algebra representations of Braid groups and link polynomials." \emph{Annals of Mathematics} \textbf{126} (1987): 103-111.
\bibitem[Jo3]{Jo3} Jones, V.  ``On knot invariants related to some statistical mechanical models."  \emph{Pacific Journal of Mathematics} \textbf{137} (1989): 311-334.
\bibitem[Ju]{Ju} Juyuyama, J. ``Sur les nouveaux g\'en\'erateur de l'alg\'ebre de Hecke $\H(G,U,1)$." \emph{Journal of Algebra} \textbf{204} (1998): 49-68.
\bibitem[KL]{KL} Kazhdan, D.  and Lusztig G.  ``Representations of Coxeter groups and Hecke algebras." \emph{Inventiones Mathematicae} \textbf{53} (1979): 165-184.
\bibitem[Kn]{Kn} Knuth, D.  ``Permutations, matrices, and generalized Young tableaux." \emph{Pacific Journal of Mathematics}  \textbf{34} (1970): 709-727.
\bibitem[LV]{LV} Lusztig, G. and Vogan, D.   ``Singularities of closures of $K$-orbits on flag manifolds." \emph{Inventiones Mathematicae}  \textbf{21} (1983): 365-379.
\bibitem[Ma]{Mac} Macdonald, I.G.  \emph{Symmetric Functions and Hall Polynomials,  2nd edition.}  Oxford: Oxford Science Publications, 1995.
\bibitem[St]{St} Steinberg, R.  \emph{Lectures on Chevalley Groups}.  mimeographed notes,  Yale University, 1967.
\bibitem[Yo1]{Y1} Yokonuma, T.  ``Sur le commutant d'une repr\'esentation d'un groupe de Chevalley fini." \emph{Journal of the Faculty of Science, University of Tokyo}  \textbf{15}  (1968): 115-129.
\bibitem[Yo2]{Y2} Yokonuma, T.  ``Sur le commutant d'une repr\'esentation d'un groupe de Chevalley fini II."  \emph{Journal of the Faculty of Science, University of Tokyo}  \textbf{16} (1969): 65-81.
\bibitem[Ze]{Z} Zelevinsky, A.  \emph{Representations of Finite Classical Groups}.  New York: Springer Verlag, 1981.
\end{thebibliography}
\end{document}